\documentclass[letterpaper, 10 pt, conference]{ieeeconf}  

\IEEEoverridecommandlockouts                              

\usepackage{cite}
\usepackage{amsmath,amssymb,amsfonts}
\usepackage{algorithm2e}
\usepackage{graphicx}
\usepackage{textcomp}
\usepackage{xcolor}
\usepackage{caption}    
\usepackage{subcaption} 
\newcommand{\tr}{\mathsf{T}}  
\newcommand{\norm}[1]{\left\lVert#1\right\rVert}
\usepackage{url} 
\usepackage{bm} 
\usepackage{placeins} 
\newcommand{\Tr}{\mathrm{Tr}} 
\DeclareMathOperator*{\argmin}{arg\,min}
\def\BibTeX{{\rm B\kern-.05em{\sc i\kern-.025em b}\kern-.08em
    T\kern-.1667em\lower.7ex\hbox{E}\kern-.125emX}}
\begin{document}

\title{Generalized Maximum Entropy Differential Dynamic Programming \\
 }

\author{Yuichiro Aoyama$^{1,2}$ and Evangelos A. Theodorou$^{1}$
\thanks{$^{1}$School of Aerospace Engineering, Georgia Institute of Technology, Atlanta, GA, USA}
\thanks{$^{2}$Komatsu Ltd.,Tokyo, Japan}
\thanks{\{yaoyama3, evangelos.theodorou\}@gatech.edu}}
\maketitle

\begin{abstract}
We present a sampling-based trajectory optimization method derived from the maximum entropy formulation of Differential Dynamic Programming with Tsallis entropy. This method is a generalization of the legacy work with Shannon entropy, which leads to a Gaussian optimal control policy for exploration during optimization. With the Tsallis entropy, the policy takes the form of $q$-Gaussian, which further encourages exploration with its heavy-tailed shape. Moreover, the sampling variance is scaled according to the value function of the trajectory. This scaling mechanism is the unique property of the algorithm with Tsallis entropy in contrast to the original formulation with Shannon entropy, which scales variance with a fixed temperature parameter. Due to this property, our proposed algorithms can promote exploration when necessary, that is, the cost of the trajectory is high. The simulation results with two robotic systems with multimodal cost demonstrate the properties of the proposed algorithm.

\end{abstract}

\section{Introduction}
Tsallis entropy, also known as $q$-logarithmic entropy, is a generalization of the standard Boltzmann-Gibbs or Shannon entropy \cite{Tsallis1988Tsallisentropy,{Leinster2021qentropy}}.
The Tsallis entropy is nonadditive, which means that the sum of the entropy of probabilistically independent subsystems is not equal to that of the entire system. The entropy is used in nonextensive statistical mechanics which can handle strongly correlated random variables\cite{Tsallis2005extensive, Umarov2008nonextensive}. It also helps to analyze complex phenomena in physics such as in \cite{Lutz2003opticallattice, Liu2008Super}, etc. 

Maximization of Shannon entropy under fixed mean and covariance yields a Gaussian distribution \cite{bishop2006pattern}. 
With the Tsallis entropy, similar maximization with a fixed $q$-mean and covariance leads to a $q$-Gaussian distribution \cite{Tsallis1995, Prato1999Nonextensive}.  The distribution has a heavy tail compared to a normal Gaussian with a certain range of the entropic index $q$. Due to this property, the $q$-Gaussian distribution is utilized for several engineering applications. In \cite{INOUE2013qgaussimage}, a mixture of the distribution is used for image and video semantic mapping, improving the robustness to outliers. In stochastic optimization, the distribution is used as a generalization of the Gaussian kernel, showing better control ability in smoothing functions\cite{Ghoshodastidar2012qGauss}. It was also used for mutation in the evolutionary algorithm, showing its effectiveness over other distributions \cite{Tinos2011qgaussevolution}. 

Maximum Entropy (ME) is a popular technique in various fields, such as Stochastic Optimal Control (SOC) and Reinforcement Learning (RL), which can improve the robustness of stochastic policies. In ME, an entropic regularization term  is added to the objective, encouraging exploration during optimization and preventing policies from converging to a delta distribution.\cite{Haarnoja2017deepenergy, Ziebart2010maxentropy}. This is because the regularized objective simultaneously maximizes entropy while minimizing the original objective.
In SOC, Shannon entropy, which leads to Kullback-Leibler (KL) divergence, is used as a regularization term between the controlled and prior distributions \cite{Williams2017MPPI}. In \cite{wang2021variational}, Tsallis divergence, a generalization of the KL divergence, is used, showing improvements in robustness. For RL application, Tsallis entropic regularization leads to better performance and faster convergence \cite{chen2018effectiveTsallisRL, Lee2019genentropyRL, zhu2022enforcingTsallisRL}. Recently, the ME technique with Shannon entropy has been applied to a trajectory optimization algorithm Differential Dynamic Programming (DDP) \cite{Jacobson1970ddp}, which yields a new algorithm ME-DDP \cite{so2022meddp}. DDP is a powerful
trajectory optimization tool that has a quadratic convergence rate \cite{liao1991convergence}. However, since it relies on local information, of cost and dynamics, it converges to local minima rather than global minima. In contrast, ME-DDP can explore multiple local minima while minimizing the original objective. 
Consequently, it can find better local minima than those of normal DDP. 

In this paper, we propose ME-DDP with Tsallis entropy, which is a generalization of ME-DDP with Shannon entropy. The Tsallis entropy turns the optimal policy from Gaussian to heavy-tailed $q$-Gaussian, improving the exploration capacity of ME-DDP. Moreover, the generalized formulation can scale the variance of the policy based on the value function of the trajectory, which further promotes exploration. We validate our proposed algorithm in two robotic systems, i.e., a 2D car and a quadrotor, and make comparison with normal DDP and ME-DDP with Shannon entropy.

The main contribution of this paper is as follows.
\begin{itemize}
\item We derive DDP with Tsallis entropic regularization.
\item We show the superior exploration capability of ME-DDP with Tsallis entropy to ME-DDP with Shannon entropy by analyzing the stochastic policies.
\item We validate the exploration capability of our method with two robotic systems in simulation.
\end{itemize}



\section{Preliminaries}
\subsection{Tsallis Entropy}
With an entropic index $q \in \mathbb{R}$, we introduce $q$-logarithm function 
\begin{align}\label{eq:def_q_log}
\ln_{q}(x) =
\begin{cases}
\ln (x), \quad \quad \quad &q =1, \ x>0\\
\frac{x^{1-q}-1}{1-q}, \quad  &q \neq 1 \ x>0, 
\end{cases}
\end{align}
and its inverse $q$-exponential function 
\begin{align*}
\exp_{q}(x) = 
\begin{cases}
\exp(x), \hspace{30mm} &q=1,\\ 
[1-(q-1)x]_{+}^{-\frac{1}{q-1}},  &q\neq 1. 
\end{cases}
\end{align*}
Consider a discrete set of probabilities $\{p_i\}$. $i = 1, \cdots, I$. The Tsallis or $q$-entropy is defined as 
\begin{align}\label{eq:tsallis_entropy_def}
S_{q}(p) = \frac{1-\sum_{i}p_{i}^{q}}{q-1},
\end{align}
which is a generalization of the Shannon entropy
\begin{align*}
H(p) = \sum_{i}p_{i}\log\Big(\frac{1}{p_{i}}\Big) = -\sum_{i}p_{i}\log{p_{i}}.
\end{align*}
The Tsallis entropy may be represented as a sum of the product of probability and $q$-$\log$ probability as 
\begin{align*}
 S_{q}(p) &= \sum_{i}p_{i}\ln_{q}\Big(\frac{1}{p_{i}}\Big)
= \sum_{i}p_{i}\frac{p_{i}^{q-1}-1}{1-q}
= \frac{1-\sum_{i}p_{i}^{q}}{q-1},
\end{align*}
which recovers the definition in \eqref{eq:tsallis_entropy_def}. We believe that this is the most straightforward generalization. Although 
$$S_{q}(p) = -\sum_{i}p_{i}\ln_{q}p_{i}$$ 
can be another option for the definition, the following relation 
\begin{align*}
-\ln_{q} p_{i} \neq \ln_{q}(1/p_{i}),
\end{align*}
changes the parametrization from $q$ to $2-q$ as in \cite{Lee2019genentropyRL}.
We use the former definition in \eqref{eq:tsallis_entropy_def} to avoid this change. The argument above holds with continuous probability distribution by changing the summation to integral.

\subsection{Univariate q-Gaussian distribution}
Let us consider a Probability Distribution Function (PDF) $p_{q}(x)$ of a scalar random variable $x \in \mathbb{R}$. $q$-Gaussian distribution is a generalized version of Gaussian distribution, which is obtained by maximizing the $q$-entropy with a fixed (given) $q$-mean $\mu_{q}$ and $q$-variance $\sigma^{2}_{q}$, where the subscript $q$ means that they are computed with the $q$-escort distribution, which is a normalized $q$ th power of the original $p_{q}(x)$, i.e.,
\begin{align*}
\mu_{q} = \frac{\int xp_{q}(x)^{q} \mathrm{d}x}{\int p_{q}(x)^{q} \mathrm{d}x}, \quad 
\sigma_{q}^{2} = \frac{\int (x-\mu_{q})^{2}p_{q}(x)^{q} \mathrm{d}x}{\int p_{q}(x)^{q} \mathrm{d}x}.
\end{align*}
This normalization is known to be the correct formulation of nonextensive statistics\cite{Prato1999Nonextensive}. 
The univariate $q$-Gaussian distribution is given as follows \cite{Prato1999Nonextensive}.
\begin{align*}
p_{q}(x)&=\frac{1}{Z_{q}} \Big(1 - \frac{1-q}{3-q}\frac{(x-\mu_{q})^{2}}{\sigma_{q}^{2}}\Big)^{\frac{1}{1-q}}, \\
\text{with} \
Z_{q} &= 
\begin{cases}
\big[\sigma_{q}^{2}\frac{3-q}{q-1}\big]^{\frac{1}{2}} B\big( \frac{1}{2}, \frac{3-q}{2(q-1)}\big), \ 1<q<3,\\
\big[\sigma_{q}^{2}\frac{3-q}{1-q}\big]^{\frac{1}{2}} B\big( \frac{1}{2}, \frac{2-q}{1-q}\big), \ q<1,
\end{cases}
\end{align*}
where $B(\cdot, \cdot)$ is Beta function.
The bounds for $q$ are for the convergence of integrals to compute the partition function. Here, we provide two properties of the function.\\
\textbf{Compact support when $q<1$.}\\
To satisfy $p_{q}(x)\geq 0$, $p_{q}(x)$ has a compact support in the case of $q<1$, which is given by $
x \in [x_{b}, x_{u}]$, with
\begin{align*}
x_{b} = \mu - a_{q}, \
x_{u} = \mu + a_{q}, \
a_q = \sqrt{{(3-q)}/{(1-q)}}.
\end{align*}
Note that $\mu_{q}=\mu$ when $q<2$.\\
\textbf{Recovers Normal Gaussian as $q \rightarrow 1$.}\\
Using the definition of exponential
$e^{x} = \lim_{n \rightarrow 0} (1+x/n)^{n}$,
$p_{q}$ approaches 
$$p_{1}(x) = (1/\sqrt{2\pi}\sigma) \exp[(x-\mu)^{2}/{2\sigma^{2}}],$$
as $q \rightarrow 1$, which is normal Gaussian distribution. In addition to these, moments are only defined with certain $q$s due to the convergence condition of integrals. A sampling method from the distribution is proposed in \cite{Thistleton2007genBoxmuller}.
\begin{figure}[t]
\centerline{{\includegraphics[trim={0.6cm 0.5cm 0cm 0cm},clip,width=\linewidth]{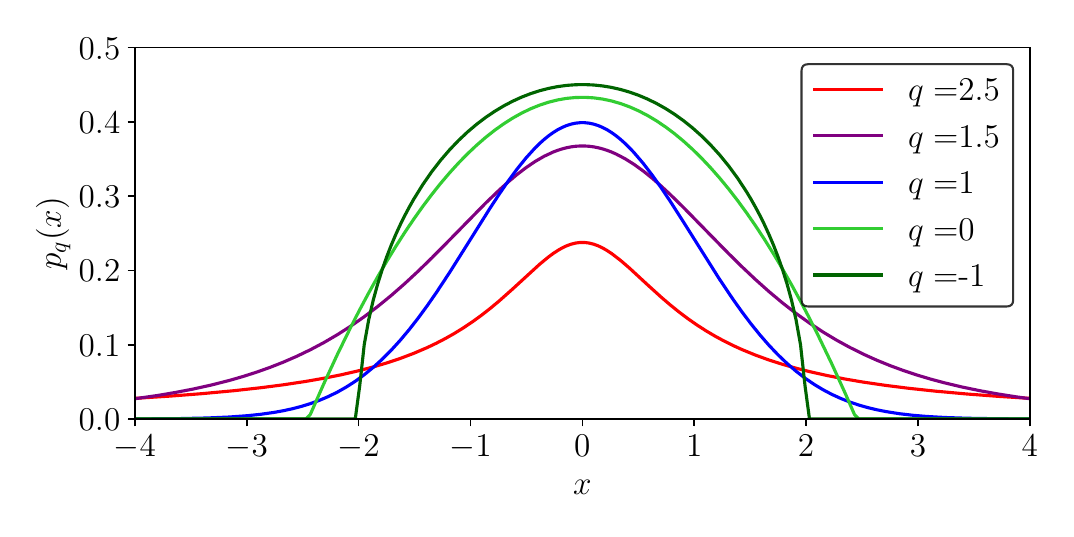}}}
\caption{$q$-Gaussian distribution with different $q$s with $\mu_{q}=0$ and $\sigma_{q}^{2}=1$. $q=1$ corresponds to a normal Gaussian.}
\label{fig:q_distribution}
\end{figure}

We show some $p_{q}(x)$ with different $q$ in Fig. \ref{fig:q_distribution}. As can be seen from the figure, a large $q$ with $q>1$ yields a heavier-tailed distribution than the normal Gaussian ($q=1$). We can also observe the compact support when $q<1$. For exploration purposes, we prefer the heavy-tailed distribution and thus focus on the case where $q>1$, hereafter.
\subsection{Multivariate q-Gaussian Distribution}
The multivariate variant of $p_{q}$ for $q>1$ and $x\in\mathbb{R}^{n}$ is
\begin{align}\label{eq:multi_qgauss}
\notag
&p_{q}(x) = \frac{1}{Z_{q}} \Big[1+\frac{q-1}{n+2-nq}(x-\mu_q)^{\tr}\Sigma_{q}^{-1}(x-\mu_q)]^{-\frac{1}{q-1}},\\
&\text{with} \
Z_{q} = \Big[ \frac{(n+2)-nq}{q-1}\Big]^{\frac{n}{2}}|\Sigma_{q}|^{\frac{1}{2}}\pi^{\frac{n}{2}}\frac{\Gamma \big(\frac{1}{q-1}-\frac{n}{2}\big)}{\Gamma\big(
\frac{1}{q-1}\big)},
\end{align} 
where $\Gamma(\cdot)$ is Gamma function \cite{VIGNAT2006167,Vignat2007deformedexp}, and $\Sigma_{q} \in \mathbb{R}^{n\times n}$ is $q$-covariance, which is a multivariate version of $\sigma_{q}^{2}$. Here, $q$ needs to satisfy the following condition. 
\begin{align}\label{eq:requirement_q}
1 < q < 1 + {2}/{n}.
\end{align}
In addition, when the tighter condition 
\begin{align}\label{eq:requirement_q_multi_cov}
1 < q < (n+4)/(n+2) = 1 + 2/(n+2),
\end{align}
is satisfied, (normal, not $q$) covariance $\Sigma \in \mathbb{R}^{n\times n}$ exists and is finite. It has the following relation with $\Sigma_{q}$   
:
\begin{align*}
\Sigma = \frac{n+2-nq}{n+4-(n+2)q}\Sigma_{q},
\end{align*}
whose coefficient is positive due to \eqref{eq:requirement_q}, \eqref{eq:requirement_q_multi_cov}. Under $1<q<1  +2/(n+1)$, which is slightly looser than \eqref{eq:requirement_q_multi_cov}, $\mu = \mu_{q}$.
Eq. \eqref{eq:multi_qgauss} corresponds to the Student's $t$ distribution\cite{Student1908tdist, bishop2006pattern}
\begin{align*}
&p_{t}(x) \\
=& \frac{\Gamma(\nu+n)/2}{\Gamma(\nu/2)\nu^{n/2}\pi^{n/2}|\Sigma_{t}|} \Big[1+\frac{1}{\nu}(x-\mu_{t})^{\tr}\Sigma_{t}^{-1}(x-\mu_{t})\Big]^{-\frac{\nu+n}{2}}.
\end{align*}
$\nu$ here is known as degrees of freedom. By taking 
\begin{align}\label{eq:scale_student_q_gauss}
\nu = \frac{n+2-nq}{q-1}, \quad \frac{q-1}{n+2-nq} \Sigma_{q}^{-1} = \frac{1}{\nu} \Sigma_{t}^{-1},
\end{align}
the $q$-Gaussian distribution $(q>1)$ is recovered. 
We use this property to sample from the distribution in section \ref{sec:numerical_exp}.

\subsection{Maximum Entropy DDP}
In this section, we have a brief review of ME-DDP with unimodal and multimodal policies \cite{so2022meddp}. Consider a trajectory optimization problem of a dynamical system with state $x\in \mathbb{R}^{n_x}$ and control $u\in \mathbb{R}^{n_u}$. Let us define the state and control trajectory with time horizon $T$ as
\begin{align*}
X = [x_{0}, \cdots, x_{T}], \quad U = [u_{0}, \cdots, u_{T-1}],
\end{align*}
and deterministic dynamics $x_{t+1} = f(x_{t}, u_{t})$.
The problem is formulated as a minimization problem of the cost  
\begin{align}\label{eq:DDP_cost}
J(X,U) = J(x_{0}, U) = \sum_{t=0}^{T-1}l_{t}(x_{t},u_{t}) + \Phi(x_{T}),
\end{align}
subject to the dynamics. Here, $l_{t}$ and $\Phi$ are running and terminal cost, respectively.
We consider a stochastic control policy $\pi(x_{t}|u_{t})$ with the same dynamics. Moreover, we add the Shannon entropic regularization term
\begin{align*}
H[\pi_{t}] = \mathbb{E}_{\pi_{t}}[\ln{\pi_{t}}]= - \int \pi_{t}(u_{t})[\ln{\pi_{t}}]\mathrm{d}u_{t},
\end{align*}
to the cost and consider the expectation of the cost over $\pi$
\begin{align*}
J_{\pi} = J_{\pi}(x_{0}, \pi) = \mathbb{E}\big[\Phi(x_{T}) + \sum_{t=0}^{T-1}\big(l_{t}(x_{t},u_{t}) - \alpha H[\pi_{t}]\big)\big],
\end{align*}
where $\alpha(>0)$ is a temperature that determines the effect of the regularization term \cite{Haarnoja2017deepenergy}. 
Applying the Bellman's principle for normal DDP, i.e., 
\begin{align*}
V(x_{t}) = \min_{u_{t}} \{ l_{t}(x_{t}, u_{t}) + V(x_{t+1})\},
\end{align*}
with value function $V$, in our setting, we have
\begin{align}\label{eq:MEDDP_V_definition}
V(x_{t}) = \min_{\pi_{t}}\big\{ \mathbb{E}_{\pi}[l_{t}(x_{t}, u_{t}) + V(x_{t+1})] - \alpha  H[\pi_{t}] \big\}
\end{align}
The solution to right-hand side of the equation is
\begin{align}\label{eq:max_entropy_with_cnst}
\notag
\min_{\pi} \mathbb{E}_{u \sim \pi(\cdot|x)}[\underbrace{l(u,x) + V(x')}_{Q(x,u)}] - \alpha H[\pi(\cdot | x)], \\ 
\text{subject to} 
\int \pi(u|x){\mathrm{d}}u = 1,
\end{align}
where we dropped time instance $t$ and denote $x_{t+1}$ as $x'$. $\pi(\cdot | x)$ means that it originally was a function of $u$, but it vanishes after taking the expectation to compute $H$. We emphasize that the expectation is computed by sampling the controls from $\pi$ by changing the notation in the expectation. 
The optimization problem above is solved by forming the Lagrangian and setting the functional derivative of $\pi$ zero. The optimal policy and value function is obtained as
\begin{align}
\pi^{\ast}(u|x) &= \frac{1}{Z(x)}\exp\big[ -\frac{1}{\alpha}Q(x,u)\big], \label{eq:ME_DDP_pi_star }\\
V(x) &= -\alpha \ln Z(x) \label{eq:MEDDP_V=-lnZ},
\end{align}
with a partition function $Z(x)$.
\subsubsection{Unimodal policy}
To obtain $\pi^{\ast}$,
we first consider a deviation $\delta x, \delta u$ from a nominal trajectory $\bar{x}, \bar{u}$, having a pair $x = \bar{x} + \delta x, u = \bar{u} + \delta u$. Then, we perform a quadratic approximation of $Q$ around $\bar{x}$, $\bar{u}$, and plug it in \eqref{eq:ME_DDP_pi_star }, having a Gaussian policy 
\begin{align}\label{eq:unimodal_pi_gaussian}
\pi^{\ast}(\delta u|\delta x) \sim \mathcal{N}(\delta u^{\ast}, \alpha Q_{uu}^{-1}),
\end{align}
where $\delta u^{\ast}$ is a solution of normal DDP, i.e.,
\begin{align}\label{eq:delta-u-star}
     \delta {u}^{\ast}&= k+K\delta x, \\\notag
    \text{with}\ k &= -{Q}^{-1}_{{uu}}{Q_{{u}}},\ K = -{Q}^{-1}_{{uu}}{Q_{{ux}}}.
\end{align}

As in normal DDP, ME-DDP has backward and forward passes. It can have $N$ (two in the original work) trajectories in parallel.
The backward pass is the same as that of the normal DDP. In the forward pass, a new control sequence is sampled from the optimal policy based on the best (lowest cost) trajectory for every $m$ iteration. In the sampling phase, the best trajectory is kept and only the remaining $N-1$ trajectories are sampled. Aside from that, the pass is the same as normal DDP, that is, forward propagation of dynamics and line search for cost reduction.

\subsubsection{Multimodal policy}\label{subsec:MEDDP_mul_policy}
Here, we consider LogSumExp approximation of the value function using $N$ trajectories, which gives the terminal state of the value function as
\begin{align*}
V(x_{T}) = \tilde{\Phi}(x)
= -\alpha \ln \sum_{n=1}^{N} \exp\Big[
-\frac{1}{\alpha}\Phi^{(n)}(x)
\Big],
\end{align*}
where $\Phi^{(n)}(x)$ $n=1, \cdots N$, is the terminal cost of $n$th trajectory.
Exponential transformation
$
\mathcal{E}_{\alpha}(y) = \exp(
-y/\alpha),
$
of the value function allows us to write $Z$ as a sum of those of $N$ trajectories denoted by $z^{(n)}$, 
\begin{align*}
z(x) &= \mathcal{E}_{\alpha}[V(x)] = \exp \big(
-(-\alpha \ln Z(x))/{\alpha}
\big)\\
&= Z(x) = \sum_{n=1}^{N} z^{(n)}(x).
\end{align*}
Let us also transform the running cost, having the desirability function $r(x,u)=\mathcal{E}_{\alpha}[l(x,u)]$.
Due to the property of the exponential, 
i.e., $\exp(x+y) = \exp(x)\exp(y)$, the optimal policy becomes
\begin{align*}
\pi(u|x) = \frac{1}{Z}\exp \Big[-\frac{1}{\alpha} (\underbrace{l(u,x) + V(x')}_{Q(u,x)})\Big] = \frac{1}{Z}z(x')r(x,u).
\end{align*}
Thanks to these properties, we obtain the optimal policy as a weighted sum of those of each trajectory $\pi^{(n)}$s as 
\begin{align*}
 \pi(u|x)
&= \sum_{n=1}^{N} w^{(n)}(x) \pi^{(n)}(u|x),
\end{align*}
with weight $w^{(n)} \propto \mathcal{E}_{\alpha}[V^{(n)}(x)]$.
Since $\pi^{(n)}$ are Gaussian as in \eqref{eq:unimodal_pi_gaussian}, the policy is now a mixture of Gaussians (and thus multimodal) with weights based on value functions.

\section{Generalized Max Entropy DDP}\label{sec:genMEDDP}
\subsection{Tsallis entropic regularization}
Based on ME-DDP in the previous section, we now consider entropic regularization with Tsallis entropy
\begin{align*}
S_{q} = \frac{1-\int p^{q}(x)\mathrm{d}x}{q-1} .
\end{align*}
Let us revisit the optimization problem in \eqref{eq:max_entropy_with_cnst}. We now use $S_{q}$ instead of $H$ as a regularization term, and consider the $q$-escort distribution of $\pi$ with a normalization constant $C$, which transforms the problem as  
\begin{align}\label{eq:gen_min_wrt_pi}
\min_{\pi} \mathbb{E}_{\pi^{q}(\cdot|x)/C}[l(u,x) + V(x')] - \alpha S_{q}[\pi(\cdot | x)],
\end{align}
under the same constraints for a valid PDF for $\pi^{q}/C$ as in \eqref{eq:max_entropy_with_cnst}. To solve this optimization problem, we form a
Lagrangian as 
\begin{align*}
\mathcal{L} = \int Q \frac{\pi^{q}}{C} \mathrm{d}u - \alpha S_{q}[\pi]  + \lambda\Big( 1 - \int \pi \mathrm{d} \pi \Big), 
\end{align*}
with a Lagrangian multiplier $\lambda$.
The first-order optimality condition, i.e., $\nabla_{\pi} \mathcal{L} = 0$ gives 
the optimal policy as 
\begin{align}\label{eq:gen_me_ddp_pi_star}
\pi^{\ast} = Z^{-1}\Big[(q-1)\frac{Q}{C\alpha} +1\Big]^{-\frac{1}{q-1}},
\end{align}
where $Z$ is a partition function. By plugging this back into \eqref{eq:gen_min_wrt_pi} and using \eqref{eq:MEDDP_V_definition}, 
the value function is obtained as follows.
\begin{align}
&V(x)\\\notag
= &Z^{-q} \int Q\frac{\pi^{\ast q}}{C} \mathrm{d}u - \alpha \frac{1-\int \pi^{\ast q}\mathrm{d}u}{q-1}.\\\notag
= &Z^{-q+1} \frac{\alpha}{q-1} \int \underbrace{Z^{-1}
\Big[(q-1)\frac{Q}{C\alpha} + 1\Big]^{-\frac{1}{q-1}}}_{\pi^{\ast}}\mathrm{d}u-\frac{\alpha}{q-1}\\\notag
= &-\alpha \frac{Z^{1-q}-1}{1-q}
= -\alpha \ln_{q}Z \quad (\because \eqref{eq:def_q_log}).
\end{align}
Notice that this expression is obtained by changing $\log$ in \eqref{eq:MEDDP_V=-lnZ} to $q$-$\log$, which implies that the derivation is the generalization of ME-DDP from Shannon to Tsallis entropy.

\subsection{q-Gaussian Policy}\label{sec:q-gauss_policy}
To examine the property of $\pi^{\ast}$, We perform quadratic approximation of $Q$ around nominal trajectories ($\bar{x}, \bar{u}$) and complete the square as
\begin{align*}
    Q(x,u) &= {Q(\bar{x},\bar{u})} + \delta Q(\delta x, \delta u)= \tilde{V}(x) + \frac{1}{2}v^{\tr}Q_{uu}v, \\
\text{with} \quad 
{\tilde{V}(x)} &= {Q(\bar{x}, \bar{u}) + \delta Q(\delta x, \delta u^{\ast})},\
v = \delta u - \delta u^{\ast}.
\end{align*}
$\tilde{V}(x)$ is the value function of normal DDP, which is obtained by plugging $\delta u^{\ast}$ into $Q(x,u)$. $\delta Q$ contains terms up to the second order. Substituting these back into \eqref{eq:gen_me_ddp_pi_star}, we have
\begin{align}\label{eq:pi_star_q_gauss}
\pi^{\ast} &= Z^{-1}\Big[\frac{q-1}{C\alpha}\Big( \tilde{V}(x) + \frac{1}{2}v^{\tr}Q_{uu}v \Big) +1\Big]^{-\frac{1}{q-1}}\\\notag
&= Z^{-1} \Big[\frac{(q-1)\tilde{V}(x) + C\alpha}{C\alpha}\Big]^{^{-\frac{1}{q-1}}} \\\notag
\times &\Big[
1 + \frac{q-1}{n_{u}+2-n_{u}q}v^{\tr}\frac{(n_{u}+2-n_{u}q)}{2[(q-1)\tilde{V}(x) + C\alpha]}Q_{uu}v
\Big]^{-\frac{1}{q-1}},
\end{align}
which is a $q$-Gaussian (see \eqref{eq:multi_qgauss} and change of $n$ to $n_{u}$) with 
\begin{align}\label{eq:mu_q_sigma_q_q_Gausspolicy}
\mu_{q} = \delta u^{\ast}, \quad 
\Sigma_{q} = \frac{2[(q-1)\tilde{V}(x) + C\alpha]}{(n_{u}+2-n_{u}q)}Q_{uu}^{-1}. 
\end{align}
Here, $C$, the normalization term for the escort distribution, has not been determined and does not have a closed-form solution. $C$ is obtained by solving the following equation.
\begin{align}\label{eq:q_gauss_normalizer}
\int & {\pi^{q}\mathrm{d}u} ={C} \Leftrightarrow
\Big[\tilde{V}(x) + \frac{\alpha C}{q-1}\Big]^{\frac{n_{u}}{2}(q-1)}C \\\notag
& = \frac{n_{u}+2-n_{u}q}{2}\Bigg[
|Q_{uu}^{-1}|^{\frac{1}{2}}(2\pi)^{\frac{n_{u}}{2}}
\frac{\Gamma\big(\frac{1}{q-1} -\frac{n_{u}}{2}\big)}{\Gamma\big(\frac{1}{q-1}\big)}\Bigg]^{(1-q)}.
\end{align}
Since the left-hand side is monotonically increasing with $C$, and since $C$ is a scalar, the equation can be easily solved by numerical methods, such as bisection search. Observe in \eqref{eq:mu_q_sigma_q_q_Gausspolicy} that we can recover ME-DDP with Shannon entropy by $q \rightarrow 1$, which leads to $\pi^{q}/C$ ($q$-Gaussian) $\rightarrow \pi$ (Gaussian), and $\Sigma_{q} \rightarrow \alpha Q_{uu}^{-1}$.    

\begin{figure}[t]
\centerline{{\includegraphics[trim={0.6cm 0.5cm 0.45cm 0cm},clip,width=\linewidth]{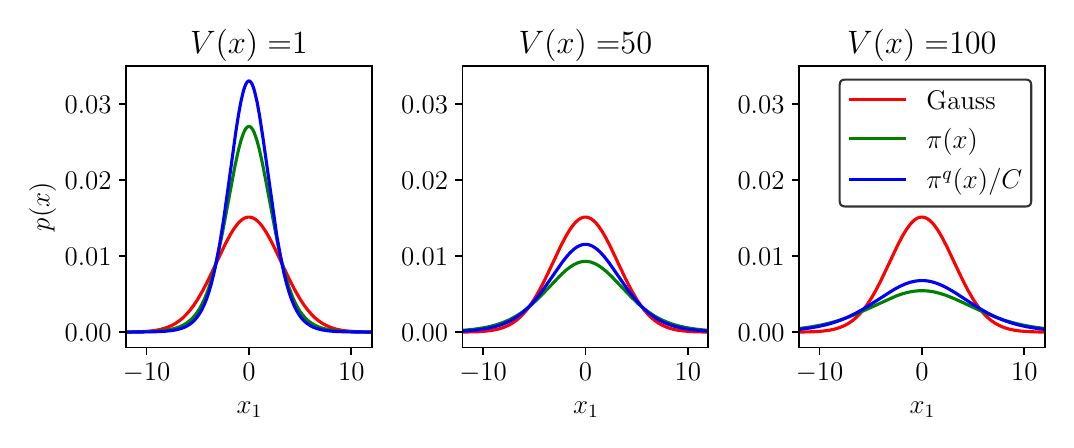}}}
\caption{Normal Gaussian policy, $q$-Gaussian policy $\pi$ and $q$-escort distribution of $\pi^{q}/C$ policy with different value function.}
\label{fig:policy_different_V}
\end{figure}

\subsection{Sampling from q-escort distribution}
In the forward pass of ME-DDP with Shannon entropy, a new control sequence is sampled from $\pi$ which is a Gaussian or a Gaussian mixture. In our case, since expectation of the cost is taken over the $q$-escort distribution of $q$-Gaussian $\pi$, it is more natural to sample from $\pi^{q}/C$ than from $\pi$. To sample from this distribution, we use the property of $q$-Gaussian, that is, the $q$-escort distribution of a $q$-Gaussian is also a $q$-Gaussian. To see this, let us introduce a parameter $q'$ and a $q$-mean and covariance as follows.
\begin{align}\label{eq:q_escort_trans_params}
q' = 2 - \frac{1}{q}, \ \mu_{q}' = \mu_{q}, \ 
\Sigma_{q}' = \frac{n_{u}+2-n_{u}q}{n_{u}+(2-n_{u})q}\Sigma_{q}.
\end{align}
With these, consider a PDF
parameterized by $q'$ as
\begin{align*}
p_{q'}(x) = \frac{1}{Z_{q'}} \Big[1-\frac{1-q'}{n_{u}+2-n_{u}q'}y^{\tr}[\Sigma_{q}']^{-1}y\Big]^{\frac{1}{1-q'}},
\end{align*}
with $y=x-\mu_{q'}$. By substituting \eqref{eq:q_escort_trans_params} in, we see that $p_{q'}(x)$ is proportional to
\begin{align*}
 \Big(\Big[1 - \frac{1-q}{n_{u}+2-n_{u}q}(x-\mu_{q})^{\tr}\Sigma_{q}^{-1}(x-\mu_{q}) \Big]^{\frac{1}{1-q}}\Big)^{q},
\end{align*}
which is a $q$ th power of $q$-Gaussian.
This implies that sampling from the $q$-escort distribution $(p_{q}(x))^{q}/C$ is achieved by sampling from the $q$-Gaussian distribution $p_{q'}(x)$ obtained by the transformation given in \eqref{eq:q_escort_trans_params}. Moreover, sampling from $q$-Gaussian is equivalent to sampling from Student's $t$ distribution with the transformation in \eqref{eq:scale_student_q_gauss} \cite{Ghoshdastidar2014smoothedfuncalgqGauss}. We use the technique to sample from the $\pi^{q}/C$ in our algorithm.

Here, we analyze the difference between the optimal control policies of ME-DDP with Shannon and Tsallis entropy. With Shannon entropy, the covariance of Gaussian $\pi^{\ast}$ is determined by $\alpha Q_{uu}^{-1}$ as in \eqref{eq:unimodal_pi_gaussian}. On the other hand, with Tsallis entropy, the covariance is also affected by the value function as in \eqref{eq:mu_q_sigma_q_q_Gausspolicy}. We can interpret the existence of the value function as follows. When the cost of the trajectory is low and thus $\tilde{V}(x)$ is low, the algorithm does not need to explore much because the current trajectory is already good. Therefore, the covariance for exploration is small. In the opposite case, the covariance is amplified by the large $\tilde{V}(x)$, which encourages exploration.
We visualize the normal Gaussian policy, $q$-Gaussian policy $\pi$, and its $q$-escort distribution $\pi^{q}/C$ with a 2D state $x = [x_{1}, x_{2}]^{\tr} \in \mathbb{R}^{2}, \alpha = 10$ and a unit $Q_{uu}$ in Fig. \ref{fig:policy_different_V} for better understanding. As shown in the top of the figure, the three panels correspond to three different $\tilde{V}(x)$ s. On the left, when $\tilde{V}(x)$ is small, $\pi^{q}$ is even tighter than Gaussian. While the Gaussian policy remains the same shape, $\pi^{q}$ becomes heavy-tailed as $\tilde{V}(x)$ increases to the right.
This means that the covarinace is properly scaled based on the cost of the trajectory, rather than scaled with the same scaling factor. Due to this property, we deduce that the ME-DDP with Tsallis entropy case has better exploration capability, which we validate in section \ref{sec:numerical_exp}.

The proposed algorithm is summarized in Alg. \ref{alg:ME_DDP}. It takes in $N$ nominal control sequences and performs optimization, outputting one pair of state control trajectories that achieve the lowest cost. In the algorithm, we write $Q_{uu}^{-1}$ as $\Sigma$ so that it can be easily compared with ME-DDP with Shannon entropy.

\RestyleAlgo{ruled}
\begin{algorithm}[hbt!]
\caption{Generalized ME-DDP}\label{alg:ME_DDP}
\SetKw{And}{and}
\KwIn{$x_{0}:$ initial state, $u^{(1:N)}$, $\Sigma^{(1:N)}$, $K^{(1:N)}$}: initial sequence, $m:$ sampling frequency, $I$: max iteration, 
$q$: entropic index that satisfies \eqref{eq:requirement_q}
.\\
\KwResult{$x^{(b)}, u^{(b)}, K^{(b)}$} 
Compute initial trajectory $(x_{0:T}^{(1:N)})$ and cost $J^{(1:N)}$.

\For{$i \gets 0$ to $I$}{
    \If{$i \% m = 0$}{
    $x^{(1)}, u^{(1)}, K^{(1)}, \Sigma^{(1)}, J^{(1)} \gets$   lowest cost mode\\    
    
   $\pi \gets$ $q$-Gaussian($u^{(1)}$, $x^{(1)}$, $K^{(1)}$, $\Sigma^{(1)}$)\;

    \For{$n \gets 2$ to $N$ in parallel}{
            $x^{(n)}, u^{(n)}, K^{(n)}, J^{(n)}$ \\
            $\gets$ Sample from $\pi^{q}/C$  
                with \eqref{eq:q_gauss_normalizer}, \eqref{eq:q_escort_trans_params}
                and sampling technique from Student's $t$ distribution in \cite{Ghoshdastidar2014smoothedfuncalgqGauss}.
                }
    }
   
    \For{$n \gets 1$ to $N$ in parallel}{
    $k^{(n)}, K^{(n)}, \Sigma^{(n)} \gets$ Backward Pass \\
    $x^{(n)}, u^{(n)}, J^{(n)} \gets $ Line Search.
    }

    }
    $b \gets \argmin_{n}J^{(n)}$
\end{algorithm}

\subsection{Availability of multimodal policy}\label{subsec:tsallis_only_uni}
In ME-DDP with Shannon entropy, the multimodal policy is available due to the additive structure of the partition function explained in section \ref{subsec:MEDDP_mul_policy}.  However, with Tsallis entropy, this is not the case. Indeed, from the partition function in \eqref{eq:pi_star_q_gauss}, we obtain
\begin{align*}
z &= \int\Big[(q-1)\frac{Q}{C\alpha} +1\Big]^{-\frac{1}{q-1}}\mathrm{d}u \\
&=\int\Big[1-(q-1)\Big(\frac{l(x,u)+V(x')}{C\alpha}\Big) \Big]^{-\frac{1}{q-1}}\mathrm{d}u \\
&= \int \exp_{q}\Big(\frac{-l(x,u)}{C\alpha} \Big) \otimes_{q} \exp_{q}\Big(\frac{-V(f(x,u))}{C\alpha}\Big)\mathrm{d}u\\
&= \int r(x,u) \otimes_{q} z(f(x,u))\mathrm{d}u,
\end{align*}
where $\otimes_{q}$ is $q$-product \cite{Vignat2007deformedexp} that is not distributed 
\begin{align*}
(a+b) \otimes_{q} c \neq
a \otimes_{q} c + b \otimes_{q} c.
\end{align*}
Therefore, $z$ is not written as a sum of $z^{(n)}$s, as opposed to the Shannon entropy case.

\begin{figure*}[!ht]
     \centering
     \begin{subfigure}[b]{\linewidth}
         \centering
         \includegraphics[trim={4cm 0.1cm 3cm 0.2cm},clip,width=\textwidth]{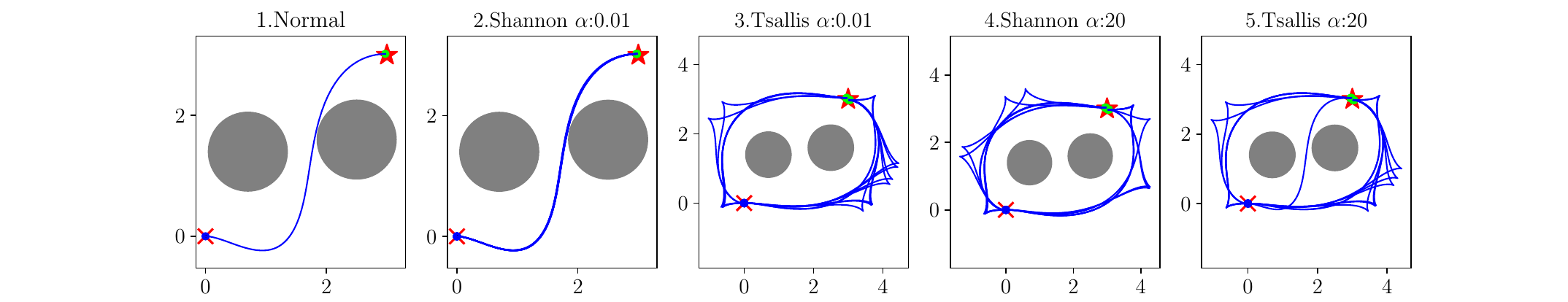}
         \caption{2D Car Example.}
         \label{fig:car_cmp}
     \end{subfigure}
     \hfill
      \begin{subfigure}[b]{\linewidth}
         \includegraphics[trim={5.3cm 1.2cm 3.3cm 0.2cm},clip,width=\textwidth]{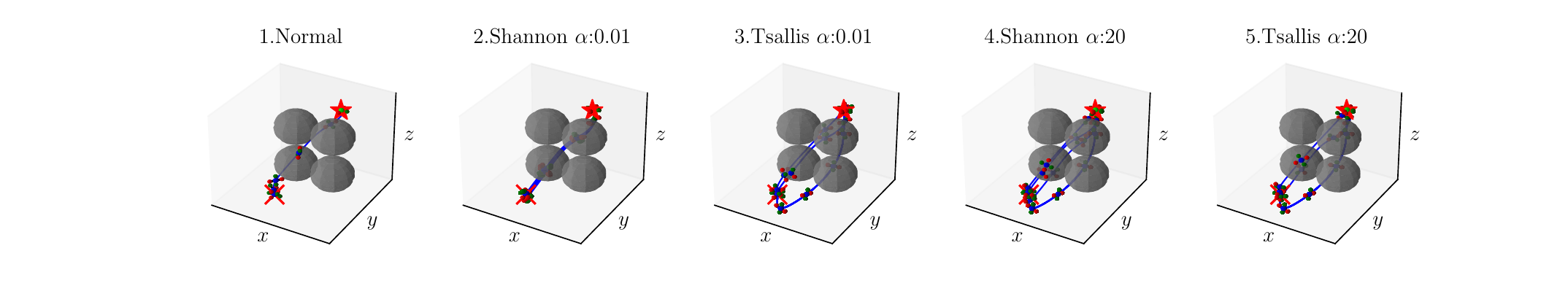}
         \caption{Quadrotor Example.}
         \label{fig:quad_cmp}
     \end{subfigure}

    \caption{Comparison of normal DDP, multimodal ME-DDP with Shannon entropy, and ME-DDP with Tsallis entropy with two different $\alpha$s. Trajectories of the 15 experiments are overlaid. {\color{red}$\times$} and {\color{red}$\bigstar$} indicate the initial and target positions, respectively. The obstacles are drawn in gray. In the quadrotor example, obstacles are drawn in transparent.} 
     \label{fig:comparison}
\end{figure*}

\begin{figure}[!ht]
     \centering
     \begin{subfigure}[b]{0.49\linewidth}
         \centering
         \includegraphics[trim={0.55cm 0cm 0.9cm 0.2cm},clip,width=\textwidth]{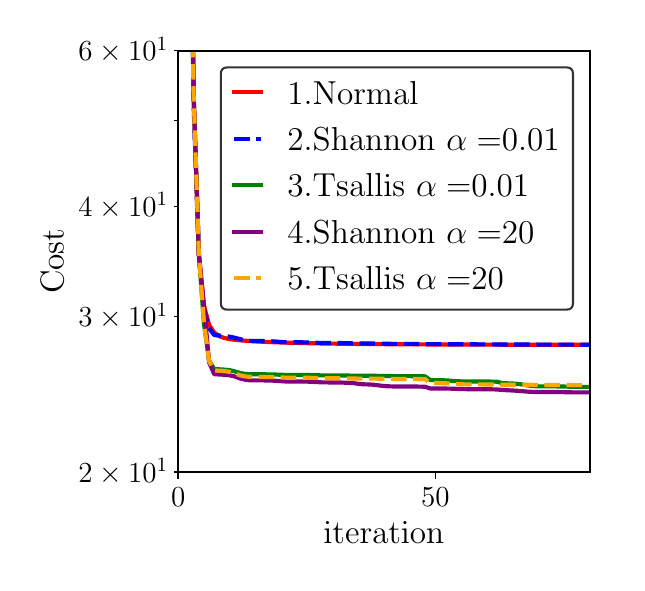}
         \caption{2D Car}
         \label{fig:car_cost}
     \end{subfigure}
     \hfill
      \begin{subfigure}[b]{0.49\linewidth}
         \includegraphics[trim={0.5cm 0cm 0.6cm 0.4cm},clip,width=\textwidth]{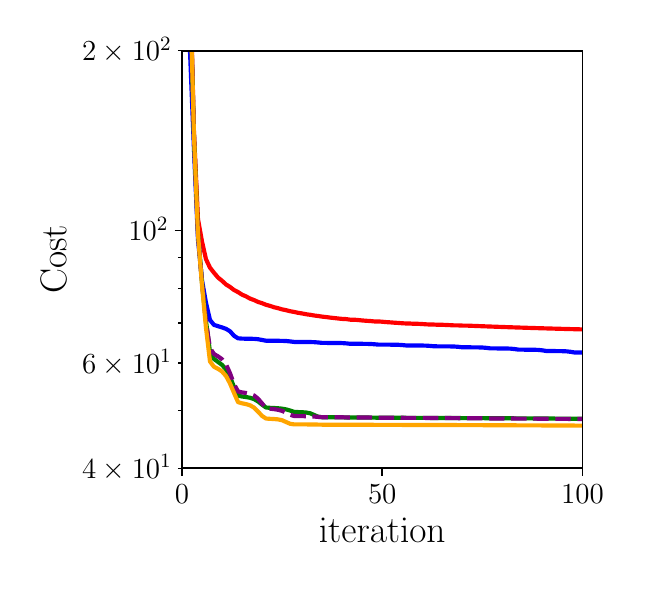}
         \caption{Quadrotor}
         \label{fig:quad_cost}
     \end{subfigure}
    \caption{Evolution of mean cost over iterations. The numbers in the legends corresponds to those in Fig.\ref{fig:comparison}}
    \label{fig:cmp_cost_change}
\end{figure}

\section{Numerical Experiments}\label{sec:numerical_exp}
In this section, we validate our proposed algorithm using two systems, a 2D car and a quadrotor. The tasks are to reach the targets while keeping distance from spherical obstacles which are encoded as part of the cost in \eqref{eq:DDP_cost}. We first give a belief description of the cost structure and systems. Then, we provide results of the experiment, including comparison with a normal DDP without entropic regularization, ME-DDP with Shannon entropy, and our proposed ME-DDP with Tsallis entropy. In the experiment, we also examine the effect of the temperature $\alpha$. Although ME-DDP with Shannon entropy has unimodal and multimodal versions, we use the multimodal one for comparison because of its superior performance \cite{so2022meddp}. We note that ME-DDP with Tsallis entropy has a unimodal policy, as explained in section \ref{subsec:tsallis_only_uni}. In both algorithms, the number of trajectories is $N=8$.

\subsection{Cost structure}
The cost for obstacles is given by 
\begin{align}\label{eq:obstacle_Gaussian}
l_{\rm{s}}(x_{t}) = \exp\big( - {(x_{t}-c_{\rm{o}})^{2}}/{2r_{\rm{o}}^{2}}\big),
\end{align}
where $c_{\rm{o}}$ and $r_{\rm{o}}$ are the center and radius of an obstacle, respectively. This cost structure is used to validate ME-DDP with Shannon entropy \cite{so2022meddp}. We use quadratic cost for the running and terminal costs in \eqref{eq:DDP_cost} and add them with \eqref{eq:obstacle_Gaussian}. The obstacle term introduces several Gaussian-shaped hills to the landscape of the cost, introducing multiple local minima.

\subsection{2D car}
The state consists of the position and orientation $\theta$ of the car, and thus the state $x = [p_{x}, p_{y},\theta]\in \mathbb{R}^{3}$. The control $u$ is translation and angular velocities $u = [v, \omega] \in \mathbb{R}^{2}$. In ME-DDP with Tsallis entropy, $q$ must satisfy $1 < q<2$ from \eqref{eq:requirement_q}. We choose $q=1.8$. The control is initialized with all zeros and thus the car stays at a initial position.

\subsection{Quadrotor}
The state of the system consists of position, velocity, orientation, and angular velocity, all of which are $\mathbb{R}^{3}$, and thus $x \in \mathbb{R}^{12}$. The control $u$ of the system is the force generated by the four rotors, which gives $u \in \mathbb{R}^{4}$. The dynamics is found in \cite{Luukkonen2011}. The requirement for $q$ is $1 <q< 1.5$, and we choose $q=1.4$. The control sequence is initialized with hovering. This sequence is obtained by setting each rotor force $mg_{0}/4$ where $m$ and $g_{0}$ are the mass of the quadrotor and gravitational acceleration.

\subsection{Results}
Fig. \ref{fig:comparison} shows the trajectories of the experiments with normal DDP, multimodal ME-DDP with Shannon entropy, and ME-DDP with Tsallis entropy with two different temperature $\alpha$s. The experiments of ME-DDPs are performed 15 times, generating 15 trajectories. Fig. \ref{fig:cmp_cost_change} shows the evolution of the mean cost of experiments over optimization iterations. 

In both systems, the normal DDP finds a local minimum that goes through between obstacles. Although it can let trajectories hit targets, the corresponding costs are high because trajectories get close to the obstacles. However, ME-DDPs can explore and find better local minima that keep a greater distance from obstacles. ME-DDP with Shanon entropy requires a large $\alpha$ to explore (see 2s and 4s in Fig. \ref{fig:comparison}). In fact, with a small $\alpha$, the results are almost the same as those of normal DDP (see 1s and 2s in Fig. \ref{fig:comparison}). This is because $\alpha$ determines the extent of exploration by scaling the sampling variance of the policy. With generalized entropy, ME-DDP can explore well even with small $\alpha$s (see 3s and 5s in Fig. \ref{fig:comparison}). Although exploration is carried out with a unimodal policy in ME-DDP with Tsallis entropy, its exploration capability is better or equivalent to the multimodal policy of ME-DDP with Shannon entropy. This is the effect of the heavy-tailed shape of $q$-Gaussian and the value function on variance \eqref{eq:mu_q_sigma_q_q_Gausspolicy} analyzed in section \ref{sec:genMEDDP}. Even with a small $\alpha$, the value function amplifies the variance for exploration when the value function (and thus the cost) of the trajectory is high, encouraging exploration. Due to this property, tuning $\alpha$ of generalized ME-DDP is easy because one can choose a small $\alpha$, rather than finding a good $\alpha$.

With a large $\alpha( = 20)$, ME-DDP with Tsallis entropy fails to find a good local minimum, which it used to find (see 3 and 5 in Fig. \ref{fig:car_cmp}). We speculate that this is because the variance is scaled too much, making it difficult to sample meaningful trajectories in the exploration process. We have observed that the ME-DDP with Shanon entropy has the same tendency with too large $\alpha$ during the experiment.


\section{Conclusion}
In this paper, we derived ME-DDP with Tsallis entropy. The algorithm has $q$-Gaussian policy whose sampling variance is scaled not only by the fixed temperature parameter $\alpha$, but also by the value function of the trajectory that is optimized. We tested our algorithm using two robotic trajectory optimization problems whose objectives have multiple local minima. The results show that although our algorithm can only have a unimodal policy, it can find better local minima even with a small temperature parameter $\alpha$ compared to the multimodal ME-DDP with Shannon entropy.  

Future work includes hardware implementation of the algorithms, as well as theoretical proofs such as the rate and condition of convergence. We are also interested in deriving ME-DDP with different entropy such as Rényi entropy.

\section*{Acknowledgement}
This work is supported by NASA under the ULI grant 80NSSC22M0070.


\begin{thebibliography}{10}
\providecommand{\url}[1]{#1}
\csname url@samestyle\endcsname
\providecommand{\newblock}{\relax}
\providecommand{\bibinfo}[2]{#2}
\providecommand{\BIBentrySTDinterwordspacing}{\spaceskip=0pt\relax}
\providecommand{\BIBentryALTinterwordstretchfactor}{4}
\providecommand{\BIBentryALTinterwordspacing}{\spaceskip=\fontdimen2\font plus
\BIBentryALTinterwordstretchfactor\fontdimen3\font minus \fontdimen4\font\relax}
\providecommand{\BIBforeignlanguage}[2]{{%
\expandafter\ifx\csname l@#1\endcsname\relax
\typeout{** WARNING: IEEEtran.bst: No hyphenation pattern has been}%
\typeout{** loaded for the language `#1'. Using the pattern for}%
\typeout{** the default language instead.}%
\else
\language=\csname l@#1\endcsname
\fi
#2}}
\providecommand{\BIBdecl}{\relax}
\BIBdecl

\bibitem{Tsallis1988Tsallisentropy}
\BIBentryALTinterwordspacing
C.~Tsallis, ``Possible generalization of boltzmann-gibbs statistics,'' \emph{Journal of Statistical Physics}, vol.~52, no.~1, pp. 479--487, Jul 1988. [Online]. Available: \url{https://doi.org/10.1007/BF01016429}
\BIBentrySTDinterwordspacing

\bibitem{Leinster2021qentropy}
T.~Leinster, \emph{Entropy and Diversity: The Axiomatic Approach}.\hskip 1em plus 0.5em minus 0.4em\relax Cambridge University Press, 2021.

\bibitem{Tsallis2005extensive}
C.~Tsallis, M.~Gell-Mann, and Y.~Sato, ``\BIBforeignlanguage{en}{Asymptotically scale-invariant occupancy of phase space makes the entropy sq extensive},'' \emph{\BIBforeignlanguage{en}{Proc Natl Acad Sci U S A}}, vol. 102, no.~43, pp. 15\,377--15\,382, Oct. 2005.

\bibitem{Umarov2008nonextensive}
\BIBentryALTinterwordspacing
S.~Umarov, C.~Tsallis, and S.~Steinberg, ``On a q-central limit theorem consistent with nonextensive statistical mechanics,'' \emph{Milan Journal of Mathematics}, vol.~76, no.~1, pp. 307--328, Dec 2008. [Online]. Available: \url{https://doi.org/10.1007/s00032-008-0087-y}
\BIBentrySTDinterwordspacing

\bibitem{Lutz2003opticallattice}
\BIBentryALTinterwordspacing
E.~Lutz, ``Anomalous diffusion and tsallis statistics in an optical lattice,'' \emph{Phys. Rev. A}, vol.~67, p. 051402, May 2003. [Online]. Available: \url{https://link.aps.org/doi/10.1103/PhysRevA.67.051402}
\BIBentrySTDinterwordspacing

\bibitem{Liu2008Super}
\BIBentryALTinterwordspacing
B.~Liu and J.~Goree, ``Superdiffusion and non-gaussian statistics in a driven-dissipative 2d dusty plasma,'' \emph{Phys. Rev. Lett.}, vol. 100, p. 055003, Feb 2008. [Online]. Available: \url{https://link.aps.org/doi/10.1103/PhysRevLett.100.055003}
\BIBentrySTDinterwordspacing

\bibitem{bishop2006pattern}
\BIBentryALTinterwordspacing
C.~Bishop, \emph{Pattern Recognition and Machine Learning}.\hskip 1em plus 0.5em minus 0.4em\relax Springer, January 2006. [Online]. Available: \url{https://www.microsoft.com/en-us/research/publication/pattern-recognition-machine-learning/}
\BIBentrySTDinterwordspacing

\bibitem{Tsallis1995}
\BIBentryALTinterwordspacing
C.~Tsallis, S.~V.~F. Levy, A.~M.~C. Souza, and R.~Maynard, ``Statistical-mechanical foundation of the ubiquity of l\'evy distributions in nature,'' \emph{Phys. Rev. Lett.}, vol.~75, pp. 3589--3593, Nov 1995. [Online]. Available: \url{https://link.aps.org/doi/10.1103/PhysRevLett.75.3589}
\BIBentrySTDinterwordspacing

\bibitem{Prato1999Nonextensive}
\BIBentryALTinterwordspacing
D.~Prato and C.~Tsallis, ``Nonextensive foundation of l\'evy distributions,'' \emph{Phys. Rev. E}, vol.~60, pp. 2398--2401, Aug 1999. [Online]. Available: \url{https://link.aps.org/doi/10.1103/PhysRevE.60.2398}
\BIBentrySTDinterwordspacing

\bibitem{INOUE2013qgaussimage}
\BIBentryALTinterwordspacing
N.~Inoue and K.~Shinoda, ``q-gaussian mixture models for image and video semantic indexing,'' \emph{Journal of Visual Communication and Image Representation}, vol.~24, no.~8, pp. 1450--1457, 2013. [Online]. Available: \url{https://www.sciencedirect.com/science/article/pii/S1047320313001855}
\BIBentrySTDinterwordspacing

\bibitem{Ghoshodastidar2012qGauss}
D.~Ghoshdastidar, A.~Dukkipati, and S.~Bhatnagar, ``q-gaussian based smoothed functional algorithms for stochastic optimization,'' in \emph{2012 IEEE International Symposium on Information Theory Proceedings}, 2012, pp. 1059--1063.

\bibitem{Tinos2011qgaussevolution}
\BIBentryALTinterwordspacing
R.~Tin{\'o}s and S.~Yang, ``Use of the q-gaussian mutation in evolutionary algorithms,'' \emph{Soft Computing}, vol.~15, no.~8, pp. 1523--1549, Aug 2011. [Online]. Available: \url{https://doi.org/10.1007/s00500-010-0686-8}
\BIBentrySTDinterwordspacing

\bibitem{Haarnoja2017deepenergy}
\BIBentryALTinterwordspacing
T.~Haarnoja, H.~Tang, P.~Abbeel, and S.~Levine, ``Reinforcement learning with deep energy-based policies,'' in \emph{Proceedings of the 34th International Conference on Machine Learning - Volume 70}, ser. ICML'17.\hskip 1em plus 0.5em minus 0.4em\relax JMLR.org, 2017, p. 1352–1361. [Online]. Available: \url{https://api.semanticscholar.org/CorpusID:11227891}
\BIBentrySTDinterwordspacing

\bibitem{Ziebart2010maxentropy}
\BIBentryALTinterwordspacing
B.~D. Ziebart, ``Modeling purposeful adaptive behavior with the principle of maximum causal entropy,'' Ph.D. dissertation, Carnegie Mellon Univ., 2010. [Online]. Available: \url{https://www.cs.cmu.edu/~bziebart/publications/thesis-bziebart.pdf}
\BIBentrySTDinterwordspacing

\bibitem{Williams2017MPPI}
\BIBentryALTinterwordspacing
G.~Williams, N.~Wagener, B.~Goldfain, P.~Drews, J.~M. Rehg, B.~Boots, and E.~A. Theodorou, ``Information theoretic mpc for model-based reinforcement learning,'' in \emph{2017 IEEE International Conference on Robotics and Automation (ICRA)}, 2017, pp. 1714--1721. [Online]. Available: \url{https://ieeexplore.ieee.org/document/7989202}
\BIBentrySTDinterwordspacing

\bibitem{wang2021variational}
\BIBentryALTinterwordspacing
Z.~Wang, O.~So, J.~Gibson, B.~Vlahov, M.~S. Gandhi, G.-H. Liu, and E.~A. Theodorou, ``Variational inference mpc using tsallis divergence,'' in \emph{Robotics Science and Systems (RSS)}, 2021. [Online]. Available: \url{https://www.roboticsproceedings.org/rss17/p073.pdf}
\BIBentrySTDinterwordspacing

\bibitem{chen2018effectiveTsallisRL}
G.~Chen, Y.~Peng, and M.~Zhang, ``Effective exploration for deep reinforcement learning via bootstrapped q-ensembles under tsallis entropy regularization,'' 2018.

\bibitem{Lee2019genentropyRL}
\BIBentryALTinterwordspacing
K.~Lee, S.~Kim, S.~Lim, S.~Choi, and S.~Oh, ``Tsallis reinforcement learning: {A} unified framework for maximum entropy reinforcement learning,'' \emph{CoRR}, vol. abs/1902.00137, 2019. [Online]. Available: \url{http://arxiv.org/abs/1902.00137}
\BIBentrySTDinterwordspacing

\bibitem{zhu2022enforcingTsallisRL}
L.~Zhu, Z.~Chen, E.~Uchibe, and T.~Matsubara, ``Enforcing kl regularization in general tsallis entropy reinforcement learning via advantage learning,'' 2022.

\bibitem{Jacobson1970ddp}
D.~H. Jacobson and D.~Q. Mayne, \emph{Differential dynamic programming}.\hskip 1em plus 0.5em minus 0.4em\relax Elsevier, 1970.

\bibitem{so2022meddp}
\BIBentryALTinterwordspacing
O.~So, Z.~Wang, and E.~A. Theodorou, ``Maximum entropy differential dynamic programming,'' in \emph{2022 International Conference on Robotics and Automation (ICRA)}.\hskip 1em plus 0.5em minus 0.4em\relax IEEE Press, 2022, p. 3422–3428. [Online]. Available: \url{https://doi.org/10.1109/ICRA46639.2022.9812228}
\BIBentrySTDinterwordspacing

\bibitem{liao1991convergence}
\BIBentryALTinterwordspacing
L.-Z. Liao and C.~Shoemaker, ``Convergence in unconstrained discrete-time differential dynamic programming,'' \emph{IEEE Transactions on Automatic Control}, vol.~36, no.~6, pp. 692--706, 1991. [Online]. Available: \url{https://ieeexplore.ieee.org/document/86943}
\BIBentrySTDinterwordspacing

\bibitem{Thistleton2007genBoxmuller}
\BIBentryALTinterwordspacing
W.~J. Thistleton, J.~A. Marsh, K.~Nelson, and C.~Tsallis, ``Generalized box–mÜller method for generating $q$-gaussian random deviates,'' \emph{IEEE Transactions on Information Theory}, vol.~53, no.~12, pp. 4805--4810, 2007. [Online]. Available: \url{https://ieeexplore.ieee.org/document/4385787}
\BIBentrySTDinterwordspacing

\bibitem{VIGNAT2006167}
\BIBentryALTinterwordspacing
C.~Vignat and A.~Plastino, ``Poincaré's observation and the origin of tsallis generalized canonical distributions,'' \emph{Physica A: Statistical Mechanics and its Applications}, vol. 365, no.~1, pp. 167--172, 2006, fundamental Problems of Modern Statistical Mechanics. [Online]. Available: \url{https://www.sciencedirect.com/science/article/pii/S0378437106000744}
\BIBentrySTDinterwordspacing

\bibitem{Vignat2007deformedexp}
\BIBentryALTinterwordspacing
------, ``Central limit theorem and deformed exponentials,'' \emph{Journal of Physics A: Mathematical and Theoretical}, vol.~40, no.~45, p. F969, oct 2007. [Online]. Available: \url{https://dx.doi.org/10.1088/1751-8113/40/45/F02}
\BIBentrySTDinterwordspacing

\bibitem{Student1908tdist}
\BIBentryALTinterwordspacing
Student, ``The probable error of a mean,'' \emph{Biometrika}, vol.~6, no.~1, pp. 1--25, 1908. [Online]. Available: \url{http://www.jstor.org/stable/2331554}
\BIBentrySTDinterwordspacing

\bibitem{Ghoshdastidar2014smoothedfuncalgqGauss}
\BIBentryALTinterwordspacing
D.~Ghoshdastidar, A.~Dukkipati, and S.~Bhatnagar, ``Smoothed functional algorithms for stochastic optimization using q-gaussian distributions,'' \emph{ACM Trans. Model. Comput. Simul.}, vol.~24, no.~3, jun 2014. [Online]. Available: \url{https://doi.org/10.1145/2628434}
\BIBentrySTDinterwordspacing

\bibitem{Luukkonen2011}
\BIBentryALTinterwordspacing
T.~Luukkonen, ``Modelling and control of quadcopter,'' \emph{Independent research project in applied mathematics, Espoo}, 2011. [Online]. Available: \url{https://sal.aalto.fi/publicaitons/pdf-files/eluu11_public.pdf}
\BIBentrySTDinterwordspacing

\bibitem{Blumenson1960spherical}
L.~E. Blumenson, ``A derivation of n-dimensional spherical coordinates,'' \emph{Am. Math. Mon.}, vol.~67, no.~1, p.~63, Jan. 1960.

\end{thebibliography}


\onecolumn
\appendix

\section{Appendix}
\renewcommand{\theequation}{\thesubsection.\arabic{equation}}
\renewcommand{\thefigure}{\thesection.\arabic{figure}} 
\renewcommand{\thetable}{\thesection.\arabic{table}}
\setcounter{figure}{0} 
\setcounter{table}{0} 
\setcounter{equation}{0}
\subsection{Mathematical Background}
This section provides a mathematical background, including beta and gamma functions and some integrals used in our work.

\subsubsection{Beta and Gamma Functions}

In the derivation of $q$-Gaussian, we use the beta and gamma functions to evaluate some integrals. For the beta function, we use the following two definitions.
\begin{subequations}
    \begin{align}
    \label{eq:beta_def_inf}
B(a,b) &= \int_{0}^{\infty}\frac{t^{a-1}}{(1+t)^{a+b}} \mathrm{d}t, \\  
\label{eq:beta_def_01} 
&= \int_{0}^{1}t^{a-1}(1-t)^{b-1} \mathrm{d}t,
     \end{align}
\end{subequations}
with $a, b >0$.
For the gamma function, we use the relation with the beta function given by 
\begin{equation}\label{eq:beta_gamma}
    B(a,b) = \frac{\Gamma{(a)}\Gamma{(b)}}{\Gamma{(a+b)}},
\end{equation}
and the following property
\begin{equation}\label{eq:gamma_prop}
\Gamma(a+1) = a \Gamma(a).
\end{equation}

\subsubsection{Useful integrals (scalar)}
Here, we derive two useful integrals in our derivation. First, we compute
\begin{align*}
I_{1} = \int_{-\infty}^{\infty}(1+x^{2})^{-m}\mathrm{d}x, \quad \text{with} \ x \in \mathbb{R}, 
\end{align*}
Using \eqref{eq:beta_def_inf}, we have 
\begin{align*}
    B\Big(\frac{1}{2}, m-\frac{1}{2}\Big)
    &= \int_{0}^{\infty}
    \frac{t^{-1/2}}{(1+t)^{m}}\mathrm{d}t \\
    & = \int_{0}^{\infty}
    \frac{1/x}{(1+x^{2})^{m}}2x\mathrm{d}x \quad (\because \text{put} \ t =x^{2}, \mathrm{d}t = 2x\mathrm{d}x)\\
    & = 2\int_{0}^{\infty}
    \frac{1}{(1+x^{2})^{m}}\mathrm{d}x\\
    &= \int_{-\infty}^{\infty}{(1+x^{2})^{-m}}\mathrm{d}x = I_{1}. \quad (\because \text{the integrated is even.})
\end{align*}
Now, with \eqref{eq:beta_gamma}, we get 
\begin{align}\label{eq:int_useful_1}
I_{1} = \int_{-\infty}^{\infty}(1+x^{2})^{-m}\mathrm{d}x = B\Big(\frac{1}{2}, m-\frac{1}{2}\Big) = \frac{\Gamma\big(\frac{1}{2}\big)\Gamma\big(m-\frac{1}{2}\big)}{\Gamma{(m)}}, \ m > \frac{1}{2}.
\end{align}
Next, we compute
\begin{align*}
I_{2} = \int_{-1}^{1}(1-x^{2})^{-m}\mathrm{d}x.
\end{align*}
Here, the sign of $x^{2}$ is minus and the limit of the integral is $[-1,1]$ in contrast to $I_{1}$. 
This is because we consider Probability Distribution Function (PDF) as a part of integrand in our application. Specifically, we consider a PDF has a form of 
$$p(x) = (1-x^{2}) \ (\geq 0),$$ 
which requires $x \in [-1,1]$. Using \eqref{eq:beta_def_01}, and following the same procedure as the previous computation, we have 
\begin{align*}
    B\Big(\frac{1}{2}, -m+1 \Big) &= \int_{0}^{1}t^{-1/2}(1-t)^{-m}\mathrm{d}t \\
    &= \int_{0}^{1}(1-x^{2})^{-m}(x^{-1})2x\mathrm{d}x \\
    &= 2 \int_{0}^{1}(1-x^{2})^{-m}\mathrm{d}x \\
    &=  \int_{-1}^{1}(1-x^{2})^{-m}\mathrm{d}x.
\end{align*}
Thus, we have 
\begin{align}\label{eq:int_useful_2}
I_{2} = \int_{-1}^{1}(1-x^{2})^{-m}\mathrm{d}x = B\Big(\frac{1}{2}, -m+1 \Big) = 
\frac{\Gamma\Big(\frac{1}{2}\Big)\Gamma(-m+1)}{\Gamma{\big(-m+\frac{3}{2}\big)}}, \ m > 1
\end{align}
\subsubsection{Useful integrals (vector)}
In this section, we generalize the integrals introduced in the previous section from scalar $x$ to vector $x \in \mathbb{R}^{n}$.
We first consider a integral
\begin{equation}
I_{3}= \int_{\mathbb{R}^{n}} (1+\norm{x}^{2})^{-m} {\mathrm{d}}x.
\end{equation}
First, we transform the coordinates to $n$-dimentional spherical coordinates \cite{Blumenson1960spherical}, which has a radial coordinate $r$ and $n-1$ angular coordinates $\phi_{1}, \phi_{2}, \cdots, \phi_{n-1}$. The angles $\phi_{1}, \cdots, \phi_{n-2}$ have the range of $[0, \pi]$, and $\phi_{n}$ has $[0, 2\pi)$. With this coordinate, considering a volume element of a $n$ ball $\mathrm{d}V_{n}$, the integral is written as 
\begin{equation}
I_{3}= \int_{\mathbb{R}^{n}} (1+r^{2})^{-m} {\mathrm{d}}V_{n}.
\end{equation}
Next, we convert the volume element to the surface element. The volume element of $n$ ball is written as 
\begin{equation}
{\mathrm{d}}V_{n} = {\color{red}{r^{n-1}}}{\color{blue}\sin^{n-2}(\phi_{1})\sin^{n-3}(\phi_{2})\cdots\sin(\phi_{n-2})} {\color{red}{\mathrm{d}}r} {\color{blue}{\mathrm{d}} \phi_{1} \cdots {\mathrm{d}} \phi_{n-1}}.
\end{equation}
Similarly, we define the surface element of $(n-1)$ sphere with radius $R$ as $\mathrm{d}S_{n-1}(R)$, which is given by
\begin{equation}\label{eq:app_dS}
\mathrm{d}S_{n-1}(R) = {\color{red}R^{n-1}}{\color{blue}\sin^{n-2}(\phi_{1})\sin^{n-3}(\phi_{2})\cdots\sin(\phi_{n-2}){\mathrm{d}} \phi_{1} \cdots {\mathrm{d}} \phi_{n-1}}.
\end{equation}
From the two equations above, we can write the volume element in terms of the surface element as 
\begin{equation*}
{\mathrm{d}}V_{n} = \mathrm{d}S_{n-1}(1)r^{n-1}\mathrm{d}r.
\end{equation*}
Now, we can turn the volume integral of $I_{3}$ into the surface integral as below.
\begin{align}\label{eq:int_useful_3}
\notag
I_{3} &= \int_{\mathbb{R}^{n}}(1+r^{2})^{-m}r^{(n-1)}{\mathrm{d}}r\mathrm{d}S_{n-1} \\ \notag
& = A_{n-1}\int_{0}^{\infty}(1+s)^{-m}s^{\frac{n-1}{2}}\frac{1}{2}s^{-\frac{1}{2}}{\mathrm{d}}s \\ \notag
& = \frac{A_{n-1}}{2}\int_{0}^{\infty}(1+s)^{-m}s^{\frac{n}{2}-1}\frac{1}{2}{\mathrm{d}}s \\
& = \frac{A_{n-1}}{2}B\Big(\frac{n}{2}, m - \frac{n}{2}\Big), \quad m > \frac{n}{2} (\because \ \eqref{eq:beta_def_01})
\end{align}
where $A_{n-1}$ is the surface area of a unit sphere in $\mathbb{R}^{n}$, i.e.,
\begin{align}\label{eq:unit_sphere_are}
A_{n-1} = \frac{n \pi^{n/2}}{\Gamma\big(1 + \frac{n}{2}\big)}
= \frac{2 \pi^{n/2}}{\Gamma\big(\frac{n}{2}\big)}.
\end{align}

\renewcommand{\theequation}{\thesubsection.\arabic{equation}}
\renewcommand{\thefigure}{\thesection.\arabic{figure}} 
\renewcommand{\thetable}{\thesection.\arabic{table}}
\setcounter{figure}{0} 
\setcounter{table}{0} 
\setcounter{equation}{0}
\subsection{q-Gaussian Distribution}
This section provides how the $q$-Gaussian distribution is obtained by maximizing the entropy with fixed $q$-mean and $q$-covariance. Moreover, we derive the normal mean and covariance of the distribution. We present both univariate and multivariate cases for completeness. Although the results are available in some iterature, the derivation is not presented. Thus, we believe that the results here are worth presenting.
\subsubsection{Univariate case}
Let us consider the optimization problem defined by the following objective and constraints.
\begin{subequations}
\begin{align}
\text{Objective (Tsallis Entropy):} & \frac{1-\int p_{q}(x)^{q} \mathrm{d}x}{q-1}\\
\label{eq:const_norm}
\text{Constraint 1 (Normalization of $p_{q}(x)$):} & \int p_{q}(x) \mathrm{d}x= 1 \\
\label{eq:const_norm_q}
\text{Constraint 2 (Normalization of $p_{q}(x)^{q}$):} & \int p_{q}(x)^{q} \mathrm{d}x = C (>0)  \\
\text{Constraint 3 ($q$-mean):} & \int x p_{q}(x)^{q} \mathrm{d}x= C\mu_{q}\\
\label{eq:const_var}
\text{Constraint 4 ($q$-variance):} & \int(x-\mu_{q})^{2} p_{q}(x)^{q} \mathrm{d}x= C\sigma_{q}^{2}.
\end{align}
In this section, the integral limit is $[-\infty, \infty]$ unless otherwise specified. Fortunately, the constraint for variance includes that of the mean. Hence, the Lagrangian $\mathcal{L}$ is formed with multipliers $a, b$ as follows.
\end{subequations}
\begin{align*}
    \mathcal{L} = \frac{1-\int p(x)^{q} \mathrm{d}x}{q-1} + a\Big(\int p_{q}(x) \mathrm{d}x -1 \Big) + 
    b\Big(
    \sigma^{2}_{q}C- \int(x-\mu_{q})^{2} p_{q}(x)^{q} \mathrm{d}x \Big).
\end{align*}
Taking a functional derivative by $p_{q}(x)$, and setting it zero, we have
\begin{align*}
    \frac{-qp_{q}(x)^{q-1}}{q-1} -a -bq(x-\mu_{q})^{2}p_{q}(x)^{q-1} = 0 \\
    \Leftrightarrow
    [1-b(1-q)(x-\mu_{q})^{2}]p_{q}(x)^{q-1}=a\frac{1-q}{q}\\
    \Leftrightarrow
    p_{q}(x)^{q-1}=a\frac{1-q}{q}[1-b(1-q)(x-\mu_{q})^{2}]^{-1}.
\end{align*}
With a normalization constant $d$, we have 
\begin{align}\label{eq:q_gaussian_pq_raw}
p_{q}(x) = d[1+b(q-1)(x-\mu_{q})^{2}]^{-\frac{1}{q-1}}.
\end{align}
In order to obtain multipliers, we need to evaluate the integrals in the constraints using the integrals introduced in the previous section. Since the multiplier $b$ is for an equality constraint that is not affected by plus and minus signs, we can assume that $b>0$. We consider two cases depending on the range of $p$ to compute $b$ and $d$.\\
\noindent\underline{case 1: $q>1$.}\\
Since $b(1-q) >0$, we take 
\begin{equation}\label{eq:change_of_variable_t_b(q-1)}
    t = \sqrt{b(q-1)}(x-\mu_{q}), \quad 
    (\mathrm{d}t = \sqrt{b(q-1)} \mathrm{d} x),
\end{equation}
which from \eqref{eq:const_norm} and \eqref{eq:int_useful_1} gives 
\begin{align}\label{eq:q_gauss_const_normal}
\notag
 \int p_{q}(x) \mathrm{d}x= 1 &\Leftrightarrow 
d \int (1+t^{2})^{-\frac{1}{q-1}} \frac{1}{\sqrt{b(q-1)}}\mathrm{d}t = 1 \\ 
&\Leftrightarrow B\Big(\frac{1}{2}, \frac{3-q}{2(q-1)}\Big)[b(q-1)]^{-\frac{1}{2}} = 1,
\end{align}
Using the same change of variables for $t$, from \eqref{eq:const_norm_q}, we have
\begin{align}\label{eq:q_gauss_const_escort}
\notag
\int p_{q}(x)^{q} \mathrm{d}x= C \Leftrightarrow 
d^{q} \int (1+t^{2})^{-\frac{q}{q-1}} \frac{1}{\sqrt{b(q-1)}}\mathrm{d}t = C\\
\Leftrightarrow
B\Big(\frac{1}{2}, \frac{q+1}{2(q-1)}\Big)[b(q-1)]^{-\frac{1}{2}} = Cd^{-q}, \ (\because \eqref{eq:int_useful_1}).
\end{align}
Finally, from \eqref{eq:const_var} we have
\begin{align}\label{eq:q_gauss_const_var_t}
\notag
&\int(x-\mu_{q})^{2} p_{q}(x)^{q} \mathrm{d}x= C\sigma_{q}^{2}\\
&\Leftrightarrow
d^{q} \int \frac{t^{2}}{b(q-1)}(1+t^{2})^{-\frac{q}{q-1}} \frac{1}{\sqrt{b(q-1)}}\mathrm{d}t
= 
d^{q}[{b(q-1)}]^{-\frac{3}{2}} \int t^{2}(1+t^{2})^{-\frac{q}{q-1}}\mathrm{d}t = 
C \sigma_{q}^{2}.
\end{align}
The last integral term is computed by the change of variable $s = t^{2}$, which yields
\begin{align}\label{eq:useful_int_t2(1+t2)}
\int_{-\infty}^{\infty} t^{2}(1+t^{2})^{-\frac{q}{q-1}}\mathrm{d}t &= 2\int_{0}^{\infty} t^{2}(1+t^{2})^{-\frac{q}{q-1}}\mathrm{d}t\\\notag
&=2\int_{0}^{\infty} s(1+s)^{-\frac{q}{q-1}}\frac{1}{2}s^{-\frac{1}{2}}\mathrm{d}s \\\notag
&= \int_{0}^{\infty} s^{\frac{1}{2}} (1+s)^{-\frac{q}{q-1}}\mathrm{d}s\\\notag
&= B\Big(\frac{3}{2}, \frac{3-q}{2(q-1)}\Big). \quad (\because \eqref{eq:beta_def_inf}).
\end{align}
Note that this integral requires an additional condition $q < 3$ for convergence.
Plugging the integral back into \eqref{eq:q_gauss_const_var_t}, we get
\begin{align}
B\Big(\frac{3}{2}, \frac{3-q}{2(q-1)}\Big)[b(q-1)]^{-\frac{3}{2}} = C\sigma_{q}^{2}d^{-q}.
\end{align}
Using \eqref{eq:q_gauss_const_normal}, \eqref{eq:q_gauss_const_escort}, and \eqref{eq:q_gauss_const_var_t}, we compute multiplies and constants. Eliminating $Cd^{-q}$ using \eqref{eq:q_gauss_const_escort} and \eqref{eq:q_gauss_const_var_t} gives
\begin{align*}
    B\Big(\frac{3}{2}, \frac{3-q}{2(q-1)}\Big)[b(q-1)]^{-\frac{3}{2}}&= \sigma_{q}^{2}B\Big(\frac{1}{2}, \frac{q+1}{2(q-1)}\Big)[b(q-1)]^{-\frac{1}{2}}\\
    \Leftrightarrow
    [b(q-1)]^{-1} &= \sigma_{q}^{2}\frac{B\big(\frac{1}{2}, \frac{q+1}{2(q-1)}\big)}{B\big(\frac{3}{2}, \frac{3-q}{2(q-1)}\big)}\\
    &= \sigma_{q}^{2} \frac{\Gamma\big(\frac{1}{2}\big)
    \Gamma\big(\frac{q+1}{2(q-1)}\big)}{\Gamma\big(
    \frac{2q}{2(q-1)}\big)}  
    \frac{\Gamma \big(\frac{2q}{2(q-1)}\big)}{\Gamma\big(\frac{3}{2}\big) \Gamma\big(\frac{3-q}{2(q-1)}\big)} \quad (\because \ \eqref{eq:beta_gamma})\\
    &= \sigma_{q}^{2} \frac{\Gamma\big(\frac{1}{2}\big)
    \Gamma\big(\frac{q+1}{2(q-1)}\big)}{\Gamma\big(\frac{3}{2}\big) \Gamma\big(\frac{3-q}{2(q-1)}\big)}\\
    & = \sigma_{q}^{2} \frac{\Gamma\big(\frac{1}{2}\big)
    \frac{3-q}{2(q-1)}\Gamma\big(\frac{3-q}{2(q-1)}\big)}{\frac{1}{2}\Gamma\big(\frac{1}{2}\big) \Gamma\big(\frac{3-q}{2(q-1)}\big)} \quad (\because \eqref{eq:gamma_prop})\\
    & = \sigma_{q}^{2}\frac{3-q}{q-1} (>0, \ \because 1< q<3). 
\end{align*}
By inverting both sides, we have
\begin{align}\label{eq:b_q-1_by_sigma}
    b(q-1) = \frac{1}{\sigma_{q}^{2}}\frac{q-1}{3-q}.
\end{align}
Plugging this into \eqref{eq:q_gauss_const_normal}, the normalization constant is obtained as
\begin{align*}
d = \frac{\Big[\frac{1}{\sigma^{2}}\frac{q-1}{3-q}\Big]^{\frac{1}{2}}}{B\Big( \frac{1}{2}, \frac{3-q}{2(q-1)}\Big)}.
\end{align*}
Finally, substituting them into \eqref{eq:q_gaussian_pq_raw}, we get $p_{q}$ as
\begin{align*}
p_{q}(x)=\frac{1}{Z_{q}} \Big(1 - \frac{1-q}{3-q}\frac{(x-\mu)^{2}}{\sigma^{2}}\Big)^{\frac{1}{1-q}}, \quad 
\text{with} \quad 
Z_{q} = \Big[\sigma_{q}^{2}\frac{3-q}{q-1}\Big]^{\frac{1}{2}} B\Big( \frac{1}{2}, \frac{3-q}{2(q-1)}\Big), \ 1<q<3.
\end{align*}

\noindent\underline{case 2: $q<1$.}\\
In this case, because $b(1-q) >0$, we take 
\begin{align*}
    t = \sqrt{b(1-q)}(x-\mu_{q}),
\end{align*}
which gives the PDF as below.
\begin{equation}\label{eq:q_gaussian_pq_raw_lessthan1}
p_{q}(x) = d[1-b(1-q)(x-\mu_{q})^{2}]^{-\frac{1}{q-1}}.
\end{equation}
With this form, we evaluate the constraints. 
The constraint \eqref{eq:const_norm} gives 
\begin{align}
\notag
\int p_{q}(x) \mathrm{d}x= 1 \Leftrightarrow 
d \int (1-t^{2})^{-\frac{1}{q-1}} \frac{1}{\sqrt{b(q-1)}}\mathrm{d}t = 1\\
\label{eq:q_gauss_const_normal_lessthan1}
\Leftrightarrow d B\Big(\frac{1}{2}, \frac{q-2}{2(q-1)}\Big)[b(q-1)]^{-\frac{1}{2}} = 1, \ (\because \eqref{eq:int_useful_2} )
\end{align}
with the normalization constant $d$.
The constraint \eqref{eq:const_norm_q} yields the following result.
\begin{align}\label{eq:q_gauss_const_escort_lessthan1}
\notag
\int p_{q}(x)^{q} \mathrm{d}x= C \Leftrightarrow 
d^{q} \int_{-1}^{1} (1-t^{2})^{-\frac{q}{q-1}} \frac{1}{\sqrt{b(1-q)}}\mathrm{d}t = C \ (\because \eqref{eq:int_useful_2}) \\
\Leftrightarrow
B\Big(\frac{1}{2}, \frac{1}{1-q}\Big)[b(1-q)]^{-\frac{1}{2}} = Cd^{-q}.
\end{align}
By the similar change of variable in the previous case,
\begin{align*}
    t = \sqrt{b(1-q)}(x-\mu_{q}), \quad 
(\mathrm{d}t = \sqrt{b(1-q)} \mathrm{d} x),
\end{align*}
the constraint for variance \eqref{eq:const_var} is evaluated as follows.
\begin{align*}
d^{q} \int \frac{t^{2}}{b(1-q)}(1-t^{2})^{-\frac{q}{q-1}} \frac{1}{\sqrt{b(1-q)}}\mathrm{d}t
= 
d^{q}[{b(1-q)}]^{-\frac{3}{2}} \int t^{2}(1-t^{2})^{-\frac{q}{q-1}}\mathrm{d}t = 
C \sigma_{q}^{2}.
\end{align*}
The integral term above is computed as below:
\begin{align*}
\int t^{2}(1-t^{2})^{-\frac{q}{q-1}} \mathrm{d}t
&=
\int_{-1}^{1} t^{2}(1-t^{2})^{-\frac{q}{q-1}} \mathrm{d}t
\ \ (\because 1-t^{2} \geq 0 )
\\
&= 2 \int_{0}^{1} t^{2}(1-t^{2})^{-\frac{q}{q-1}} \mathrm{d}t\\
&= 2 \int_{0}^{1} s(1-s)^{-\frac{q}{q-1}}\frac{1}{2}s^{-\frac{1}{2}} \mathrm{d}s\\
& = \int_{0}^{1} s(1-s)^{-\frac{q}{q-1}}s^{\frac{1}{2}} \mathrm{d}s\\
&= B \Big(\frac{3}{2}, \frac{1}{1-q} \Big).
\end{align*}
We note that $\frac{1}{1-q}>0$ and thus the integral exists. By plugging this integral back, from \eqref{eq:const_var}, we get
\begin{align}\label{eq:q_gauss_const_var_lessthan1}
[{b(1-q)}]^{-\frac{3}{2}} B \Big(\frac{3}{2}, \frac{1}{1-q} \Big) = 
C d^{-q}\sigma_{q}^{2}.
\end{align}
From \eqref{eq:q_gauss_const_var_lessthan1} and \eqref{eq:q_gauss_const_escort_lessthan1}, we eliminate $d$, having
\begin{align*}
[b(1-q)] ^{-1} &= \sigma_{q}^{2}\frac{B\big(\frac{1}{2}, \frac{1}{1-q}\big)}{B\big(\frac{3}{2}, \frac{1}{1-q}\big)}\\
    &= \sigma_{q}^{2} \frac{\Gamma\big(\frac{1}{2}\big)
    \Gamma\big(\frac{1}{1-q}\big)}{\Gamma\big(
    \frac{3-q}{2(1-q)}\big)}  
    \frac{\Gamma \big(\frac{5-3q}{2(1-q)}\big)}{\Gamma\big(\frac{3}{2}\big) \Gamma\big(\frac{1}{1-q}\big)} \quad (\because \ \eqref{eq:beta_gamma})\\
    & = \sigma_{q}^{2} \frac{\Gamma\big(\frac{1}{2}\big)
    \frac{3-q}{2(1-q)}\Gamma\big(\frac{3-q}{2(1-q)}\big)}{\frac{1}{2}\Gamma\big(\frac{1}{2}\big) \Gamma\big(\frac{3-q}{2(q-1)}\big)} \quad (\because \eqref{eq:gamma_prop})\\
    & = \sigma_{q}^{2}\frac{3-q}{1-q} (>0, \because q<1). 
\end{align*}
$d$ is obtained using \eqref{eq:q_gauss_const_normal_lessthan1} as
\begin{align*}
d = \frac{\Big[\frac{1}{\sigma_{q}^{2}}\frac{1-q}{3-q}\Big]^{\frac{1}{2}}}{B\Big( \frac{1}{2}, \frac{2-q}{1-q)}\Big)}.
\end{align*}
Plugging these into \eqref{eq:q_gaussian_pq_raw_lessthan1}, $p_{q}$ is given by
\begin{align*}
p_{q}(x)=\frac{1}{Z_{q}} \Big(1 - \frac{1-q}{3-q}\frac{(x-\mu)^{2}}{\sigma_{q}^{2}}\Big)^{\frac{1}{1-q}}, \quad 
\text{with} \quad 
Z_{q} = \Big[\sigma_{q}^{2}\frac{3-q}{1-q}\Big]^{\frac{1}{2}} B\Big( \frac{1}{2}, \frac{2-q}{1-q}\Big), \ q<1.
\end{align*}
In summary, $p_{q}(x)$ is given by
\begin{align*}
p_{q}(x)=\frac{1}{Z_{q}} \Big(1 - \frac{1-q}{3-q}\frac{(x-\mu)^{2}}{\sigma_{q}^{2}}\Big)^{\frac{1}{1-q}}, \quad 
\text{with} \
Z_{q} = 
\begin{cases}
\Big[\sigma_{q}^{2}\frac{3-q}{q-1}\Big]^{\frac{1}{2}} B\Big( \frac{1}{2}, \frac{3-q}{2(q-1)}\Big), \ 1<q<3,\\
\Big[\sigma_{q}^{2}\frac{3-q}{1-q}\Big]^{\frac{1}{2}} B\Big( \frac{1}{2}, \frac{2-q}{1-q}\Big), \ q<1.
\end{cases}
\end{align*}
\textbf{Mean and Variance}\\
We derive normal (not $q$) mean and variance with $p_{q}$ obtained in the previous section.
Consider the case where $1<q$. Follwing the same change of variables as in \eqref{eq:change_of_variable_t_b(q-1)}, we have
\begin{align*}
\mu = 
\int x p_{q}(x)\mathrm{d}x &= \frac{1}{Z_{q}}\int
x[1+b(q-1)(x-\mu_{q})^{2}]^{-\frac{1}{q-1}} \mathrm{d}x \\
&= \frac{1}{Z_{q}}\int
\Big[\frac{t}{\sqrt{b(q-1)}} + \mu_{q} \Big][1+t^{2}]^{-\frac{1}{q-1}}\frac{1}{\sqrt{b(q-1)}}  \mathrm{d}t\\
&= \frac{1}{Z_{q}}\int
\frac{t(1+t^{2})^{-\frac{1}{q-1}}}{b(q-1)}\mathrm{d}t + \mu_{q}\\
& =\mu_{q}.
\end{align*}
The last equality holds when the integral in the second to the last line exists because the integrand is odd. To see when it exists, we transform the integral as follows. 
\begin{align*}
\int_{-\infty}^{\infty}
t[1+t^{2}]^{-\frac{1}{q-1}}\mathrm{d}t = 
\frac{1}{2}\int_{1}^{\infty} s^{-\frac{1}{q-1}}\mathrm{d}s.
\end{align*}
This integral exists only when $\frac{1}{q-1}>1 \Leftrightarrow q <2$. When $q$ does not meet this condition, the integral diverges, therefore, $\mu$ cannot be defined. For the variance, we have 
\begin{align*}
    \sigma^{2} &= \int(x-\mu)^{2} p_{q}(x)\mathrm{d}x\\
    &= \frac{1}{Z_{q}}\int
\frac{t^{2}}{b(q-1)} \Big( [1+t^{2}]^{-\frac{1}{q-1}}\frac{1}{\sqrt{b(q-1)}} \Big) \mathrm{d}t\\
&= \frac{1}{Z_{q}}
{b(q-1)}^{-\frac{3}{2}} \int t^{2}(1+t^{2})^{-\frac{1}{q-1}} \mathrm{d}t \\
&= \frac{1}{Z_{q}} [b(q-1)]^{-\frac{3}{2}} B\Big( \frac{3}{2}, \frac{5-3q}{2(q-1)}\Big) \ (\because \text{the same change of variable as \eqref{eq:useful_int_t2(1+t2)}})\\
&= \Big[\frac{1}{\sigma_{q}^{2}}\frac{q-1}{3-q}\Big]^{\frac{1}{2}}
\Big[{B\Big( \frac{1}{2}, \frac{3-q}{2(q-1)}\Big)}\Big]^{-1}
 \Big[\frac{1}{\sigma_{q}^{2}}\frac{q-1}{3-q}\Big]^{-\frac{2}{3}}
 B\Big( \frac{3}{2}, \frac{5-3q}{2(q-1)}\Big) \ (\because \eqref{eq:b_q-1_by_sigma})\\
 & = \frac{3-q}{q-1}{\sigma_{q}^{2}}\frac{B\big( \frac{3}{2}, \frac{5-3q}{2(q-1)}\big)}{B\big( \frac{1}{2}, \frac{3-q}{2(q-1)}\big)}\\
 & = \frac{3-q}{q-1}{\sigma_{q}^{2}} \frac{
 \frac{1}{2}\Gamma \big(\frac{1}{2}\big) 
 \Gamma\big(\frac{5-3q}{2(q-1)}\big)}{
 \Gamma \big(\frac{1}{2}\big) \frac{5-3q}{2(q-1)}\Gamma\big(\frac{5-3q}{2(q-1)}\big)
 }\\
 & = \frac{3-q}{q-1}{\sigma_{q}^{2}} \frac{q-1}{5-3q}\\
 &= \frac{3-q}{5-3q}{\sigma_{q}^{2}},
\end{align*}
which implies that it converges only when 
$$
5-3q<0 \Leftrightarrow q < \frac{5}{3}.
$$
The case with $q<1$ can be computed in a similar manner. Thus, (normal) mean and variance are given by
\begin{align}\label{eq:summary_mu_sigma_uni}
\mu = 
\begin{cases}
\mu_{q},  &q<2\\
\text{undefined}, \quad &2 \leq q < 3,
\end{cases}
\quad 
\sigma^{2} =
\begin{cases}
\frac{3-q}{5-3q}{\sigma_{q}^{2}}, & q<\frac{5}{3},\\
\infty &\frac{5}{3} \leq q < 2,\\
\text{undefined}, \quad & 2 \leq q < 3.
\end{cases}
\end{align}
Note that $\sigma$ becomes undefined when $\mu$ becomes so. This is because $\sigma$ requires $\mu$ in its definition.

\subsubsection{Multivariate Case}
In this section, we generalize the derivation in the previous section to the multivariate case with $x \in \mathbb{R}^{n}$. We consider a similar problem to the univariate case, where $q$-mean $\mu_{q}\in\mathbb{R}^{n}$ and $q$-variance is $\Sigma_{q} \in \mathbb{R}^{n \times n}$ instead of $\sigma_{q}^{2}\in\mathbb{R}$. 
\begin{subequations}
\begin{align}
\text{Objective to be minimized (Tsallis Entropy):} & \frac{1-\int p(x)^{q} \mathrm{d}x}{q-1}\\
\label{eq:const_norm_multi}
\text{Constraint 1 (Normalization of $p_{q}(x)$):} & \int p_{q}(x) \mathrm{d}x= 1 \\
\label{eq:const_norm_q_multi}
\text{Constraint 2 (Normalization of $p_{q}(x)^{q}$):} & \int p_{q}(x)^{q} \mathrm{d}x = C (>0)  \\
\text{Constraint 3 ($q$-mean):} & \int x p_{q}(x)^{q} \mathrm{d}x= C\mu_{q}\\
\label{eq:const_var_multi}
\text{Constraint 4 ($q$-covariance):} & \int(x-\mu_{q})(x-\mu_{q})^{\tr} p_{q}(x)^{q} \mathrm{d}x= C\Sigma_{q}.
\end{align}
\end{subequations}
Lagrangian of this problem is given with a multiplier $a$ and a symmetric matrix $\bm{B}$ by
\begin{align*}
\mathcal{L} = \frac{1-\int p(x)^{q} \mathrm{d}x}{q-1} + a\Big( \int p(x)\mathrm{d}x -1  \Big) +\Tr \Big(\langle \bm{B}, \Sigma_{q}-\int
(x-\mu_{q})(x-\mu_{q})^{\tr}p_{q}(x)^{q} \mathrm{d}x \rangle\Big),
\end{align*}
where $\langle \cdot,\cdot \rangle$ is a matrix inner product and $\Tr(\cdot)$ is a trace of a matrix. When a matrix constraint $\bm{A}=\bm{0} \in \mathbb{R}^{n \times n}$ exists, an inner product of the matrix and a matrix of Lagrangian multipliers $\bm{B}\in \mathbb{R}^{n \times n}$ is added it to the Lagrangian. The matrix inner product is given by $\langle \bm{B},\bm{A} \rangle = \Tr( \bm{B}^{\tr} \bm{A})$. In our setting, the constraint is given by a symmetric matrix. Thus, $\bm{B}$ is also symmetric. We use bold $\bm{B}$ to avoid confusion with the beta function. Furthermore, using the same logic as in the univariate case, we assume that $\bm{B}$ is Positive SemiDefinite (PSD). Following the same procedure as the univariate case, we take functional derivative of Lagrangian with respect to $p_{q}(x)$, and set it zero, having
\begin{align*}
\big[1 + (q-1)\Tr( \bm{B}
(x-\mu_{q})(x-\mu_{q})^{\tr}) \big]p_{q}(x)^{q-1} = a\frac{q-1}{-q}.
\end{align*}
We consider cases depending on the value of $q$.\\
\underline{case1: $q>1$}\\
The PDF has the following form:
\begin{equation}\label{eq:p_q_multi}
    p_{q}(x) = d[1+(q-1)(x-\mu_{q})^{\tr}\bm{B}(x-\mu_{q})]^{-\frac{1}{q-1}}, 
\end{equation}
where $d$ is normalization constant and the trace term is written in a quadratic form. In order to perform the integration of PDF to get $d$, we first utilize the eigenvalue and eigenvectors of $\bm{B}$, denoted by $\lambda_{i}$ and $\bm{u}_{i}$. We write an eigenvector in bold letters to emphasize that it is not the $i$ th element of a vector, but a vector itself. Since $\bm{B}$ is symmetric and real, it accepts 
eigendecomposition, which gives 
\begin{align}\label{eq:B_decomp}
\bm{B} = \sum_{i=1}^{n} \lambda_{i} \bm{u}_{i}\bm{u}_{i}^{\tr}, \quad \text{with}\ 
\bm{B}\bm{u}_{i} = \lambda_{i} \bm{u}_{i}.
\end{align}
Since $\bm{u}$ forms the basis in $\mathbb{R}^{n}$, $x-\mu_{q} (\in \mathbb{R}^{n})$ is expanded in terms of $\bm{u}$s as follows.
\begin{align}
    x - \mu = \sum_{i}^{n}y_{i}\bm{u}_{i}.
\end{align}
This expansion leads to a change of variables as
\begin{align*}
    y = \begin{bmatrix}
    y_{1} \\ \vdots \\ y_{n}
    \end{bmatrix}
    = \begin{bmatrix}
    \bm{u}_{1}^{\tr} \\
    \vdots\\
    \bm{u}_{n}^{\tr}
    \end{bmatrix}
    (x-\mu_{q}) = U(x-\mu_{q}), \quad \text{where}\quad  U = \begin{bmatrix}
    \bm{u}_{1}^{\tr} \\
    \vdots\\
    \bm{u}_{n}^{\tr}
    \end{bmatrix}.
\end{align*}
The orthonormality of $\bm{u}_{i}$ gives the following property:
\begin{equation*}
y^{\tr}U = (x-\mu_{q})^{\tr}U^{\tr}U = (x-\mu_{q})^{\tr} \Rightarrow  x-\mu_{q}  = U^{\tr}y. 
\end{equation*}
This property makes the determinant of Jacobian simple:
\begin{equation*}
\bigg\lvert \frac{\partial x}{\partial y}\bigg\rvert ^{2}= \lvert U \rvert^{2} = \lvert U \rvert \lvert U ^{\tr}\rvert = \lvert U  U^{\tr}\rvert = |I| = 1, 
\end{equation*}
and thus $\lvert {\partial x}/{\partial y}\rvert=1$.
With the new variable $y$, $d^{-1}$ is computed as
\begin{align*}
    d^{-1}&= \int_{\mathbb{R}^{n}} p_{q}(x)\mathrm{d}x \\
    &=\int_{\mathbb{R}^{n}}  [1+(q-1)(x-\mu)^{\tr}\bm{B}(x-\mu)]^{-\frac{1}{q-1}}\mathrm{d}x \\
    &= \int_{\mathbb{R}^{n}} [1+(q-1) \sum_{i=1}^{n} \lambda_{i}y_{i}^{2}]^{-\frac{1}{q-1}}\bigg\lvert \frac{\partial x}{\partial y}\bigg\rvert\mathrm{d}y \\
    &=\int_{\mathbb{R}^{n}} \Big[1+\norm{ \sqrt{q-1} \Lambda ^{\frac{1}{2}}y}^{2}\Big]^{-\frac{1}{q-1}}\mathrm{d}y,\\    
\end{align*}
where 
$
\Lambda = {\rm{diag}}{[\lambda_{1} \cdots \lambda_n]}. 
$
We then take $t = \sqrt{q-1} \Lambda ^{\frac{1}{2}}y$, whose determinant of Jacobian is given by
\begin{equation*}
\bigg\lvert \frac{\partial y}{\partial t}\bigg \rvert= \lvert (q-1)^{-\frac{1}{2}} \Lambda^{-\frac{1}{2}} \bm{I}\rvert = (q-1)^{-\frac{n}{2}} |\bm{B}|^{-\frac{1}{2}}. \quad (\because \ |\Lambda| = \prod_{i=1}^{n}\lambda_{i} = |\bm{B}|.)
\end{equation*}
This change of variable allows us to evaluate the integral of $d^{-1}$ as
\begin{align}\label{eq:q_gauss_multi_d}
    d^{-1} = \int_{\mathbb{R}^{n}}(1+\norm{t}^{2})^{-\frac{1}{q-1}}(q-1)^{-\frac{n}{2}}|\bm{B}|^{-\frac{1}{2}}\mathrm{d}t = 
    (q-1)^{-\frac{n}{2}}|\bm{B}|^{-\frac{1}{2}} \frac{A_{n-1}}{2}B\Big(\frac{n}{2}, \frac{1}{q-1}-\frac{n}{2}\Big), \ (\because \eqref{eq:int_useful_3})
\end{align}
with another condition for $q$ to make the integral converge:
\begin{align}\label{eq:multi_q_normal_cond}
\frac{1}{q-1}-\frac{n}{2}>0 \Leftrightarrow q < 1 + \frac{2}{n}.
\end{align}
By substituting \eqref{eq:p_q_multi} into \eqref{eq:const_var_multi}, we have 
\begin{align*}
C\Sigma_{q} &= d^{q}\int_{\mathbb{R}^{n}} (x-\mu_{q})(x-\mu_{q})^{\tr}[1+(q-1)(x-\mu)^{\tr}\bm{B}(x-\mu)]^{-\frac{q}{q-1}}\mathrm{d}x\\
&= d^{q} \int_{\mathbb{R}^{n}} zz^{\tr} [1+(q-1)z^{\tr}\sum_{i=1}^{n}\bm{u}_{i}\bm{u}_{i}^{\tr}z]^{-\frac{q}{q-1}}\mathrm{d}z \\
& =d^{q} \int_{\mathbb{R}^{n}} \Big(\sum_{i=1}^{n}y_{i}\bm{u}_{i}\big)
\Big(\sum_{i=1}^{n}y_{j}\bm{u}_{j}^{\tr}\Big)[1+(q-1)\sum_{k=1}^{n}\lambda_{k}y_{k}^{2}]^{-\frac{q}{q-1}}]\mathrm{d}y\\
&= d^{q}\sum_{i=1}^{n}\sum_{j=1}^{n}\int_{\mathbb{R}^{n}} \bm{u}_{i}\bm{u}_j^{\tr}[q+(q-1)\sum_{k=1}^{n}\lambda_{k}y_{k}^{2}]^{-\frac{q}{q-1}}] y_{i}y_{j} \mathrm{d}y.\\
& = d^{q}\sum_{i=1}^{n}\int_{\mathbb{R}^{n}} \bm{u}_{i}\bm{u}_i^{\tr}[1+(q-1)\sum_{k=1}^{n}\lambda_{k}y_{k}^{2}]^{-\frac{q}{q-1}}] y_{i}^{2} \mathrm{d}y
\ \  (\because \text{Orthogonality of $y_i$.})
\\
& = d^{q}(q-1)^{-\frac{n+2}{2}}|\bm{B}|^{-\frac{1}{2}}\sum_{i=1}^{n}\frac{\bm{u}_{i}\bm{u}_{i}^{\tr}}{\lambda_{i}}\int_{\mathbb{R}^{n}} (1+\norm{t}^{2})^{-\frac{q}{q-1}}t_{i}^{2}\mathrm{d}t \ \ (\because \ t=\sqrt{q-1}\Lambda^{\frac{1}{2}}y)
\\
&=d^{q}(q-1)^{-\frac{n+2}{2}}|\bm{B}|^{-\frac{1}{2}}\sum_{k=1}^{n}\frac{\bm{u}_{i}\bm{u}_{i}^{\tr}}{\lambda_{i}}I_{\Sigma}\\
& = d^{q}(q-1)^{-\frac{n+2}{2}}|\bm{B}|^{-\frac{1}{2}}I_{\Sigma} \bm{B}^{-1}, \ (\because \eqref{eq:B_decomp})
\end{align*}
where 
\begin{align*}
I_{\Sigma} &= \frac{1}{n}\int (1+\norm{t}^{2})^{-\frac{q}{q-1}}\norm{t}^{2}\mathrm{d}t \ \big(\because \int (1+\norm{t}^{2})^{-\frac{q}{q-1}}{t_{1}}^{2}\mathrm{d}t = \cdots =\int (1+\norm{t}^{2})^{-\frac{q}{q-1}}{t_{n}}^{2}\mathrm{d}t \big)\\
&= \frac{1}{n}\int_{0}^{\infty}(1+r^{2})^{-\frac{q}{q-1}}r^{2}r^{n-1} \mathrm{d}r \mathrm{d}S_{n-1} \ 
\ (\text{the same transformation as in \eqref{eq:int_useful_3}})\\
&=\frac{A_{n-1}}{n}\int_{0}^{\infty}(1+r^{2})^{-\frac{q}{q-1}}r^{n+1} \mathrm{d}r \\
&= \frac{A_{n-1}}{2n}B\Big(\frac{1}{2}n+1, \frac{1}{q-1}-\frac{n}{2}\Big). \ (\because \text{put }r = s^{2}, \eqref{eq:beta_def_01}),
\end{align*}
where the integral in the second to the last line converges due to \eqref{eq:multi_q_normal_cond}.
Thus, the constraint for variance is 
\begin{align}\label{eq:q_gauss_const_var}
d^{q}(q-1)^{-\frac{n+2}{2}}|\bm{B}|^{-\frac{1}{2}}\frac{A_{n-1}}{2n}B\Big(\frac{1}{2}n+1, \frac{1}{q-1}-\frac{n}{2}\Big)\bm{B}^{-1} = C\Sigma_{q}.
\end{align}
By substituting \eqref{eq:p_q_multi} into the constraint in \eqref{eq:const_norm_q_multi}, we get
\begin{align}\label{eq:q_gauss_const_norm_q_multi}
\notag
C &= d^{q}\int [1+(q-1)(x-\mu)^{\tr}\bm{B}(x-\mu)]^{-\frac{q}{q-1}}\mathrm{d}x\\
\notag
&= d^{q}\int [1+\norm{t}^{2}]^{-\frac{q}{q-1}}(q-1)^{-\frac{n}{2}}|\bm{B}|^{-\frac{1}{2}}\mathrm{d}t\\
\notag
&= d^{q}(q-1)^{-\frac{n}{2}}|\bm{B}|^{-\frac{1}{2}}\int [1+r^{2}]^{-\frac{q}{q-1}}r^{n-1}\mathrm{d}r \ (\because \eqref{eq:int_useful_3})\\
& = d^{q}(q-1)^{-\frac{n}{2}}|\bm{B}|^{-\frac{1}{2}}\frac{A_{n-1}}{2}B\Big(\frac{1}{2}n, \frac{q}{q-1}-\frac{n}{2}\Big),
\end{align}
where the integral in the second to the last line converges due to \eqref{eq:multi_q_normal_cond} and $q>1$.
Plugging \eqref{eq:q_gauss_const_norm_q_multi} into and \eqref{eq:q_gauss_const_var} yields
\begin{align*}
    B\Big(\frac{1}{2}n, \frac{q}{q-1}-\frac{n}{2}\Big)\Sigma_{q} &= \frac{1}{n}(q-1)^{-1}B\Big(\frac{1}{2}n+1, \frac{q}{q-1}-\frac{n+2}{2}\Big)\bm{B}^{-1}\\
\frac{1}{n}(q-1)^{-1}\bm{B}^{-1} &=\frac{B\Big(\frac{1}{2}n, \frac{q}{q-1}-\frac{n}{2}\Big)}{B\Big(\frac{1}{2}n+1, \frac{q}{q-1}-\frac{n+2}{2}\Big)}\Sigma_{q}
= \frac{\Gamma \Big(\frac{1}{2}n\Big) \Gamma \Big(\frac{q}{q-1}-\frac{n}{2}\Big)}{\Gamma\Big(\frac{1}{2}n+1\Big)\Gamma\Big(
\frac{q}{q-1}-\frac{n+2}{2}
\Big)}\Sigma_{q}
= \frac{2q-(n+2)(q-1)}{n(q-1)}\Sigma_{q},
\end{align*}
which gives the Lagrangian multiplier matrix 
\begin{align}\label{eq:q_gaussB}
\bm{B} = \frac{\Sigma_{q}^{-1}}{2q-(n+2)(q-1)} = \frac{\Sigma_{q}^{-1}}{(n+2)-nq}.
\end{align}
Plugging this back into \eqref{eq:q_gauss_multi_d} gives the normalization term
\begin{align}\label{eq:multivariate_q_d-1}
d^{-1} = \Big[ \frac{(n+2)-nq}{q-1}\Big]^{\frac{n}{2}}|\Sigma_{q}|^{\frac{1}{2}}\pi^{\frac{n}{2}}\frac{\Gamma \Big(\frac{1}{q-1}-\frac{n}{2}\Big)}{\Gamma\Big(
\frac{1}{q-1}\Big)},
\end{align}
where we use \eqref{eq:beta_gamma} to transform beta function to gamma functions.
Therefore, we obtain 
\begin{align*}
p_{q}(x) = \frac{1}{Z_{q}} \Big[1+\frac{q-1}{n+2-nq}(x-\mu_q)^{\tr}\Sigma_{q}^{-1}(x-\mu_q)\Big]^{-\frac{1}{q-1}},\quad \text{with} \ Z_{q} = d^{-1}, \ 1<q < 1+\frac{2}{n}.
\end{align*}
\underline{case2: $q<1$}\\
In this case, PDF takes the form of 
\begin{equation}\label{eq:p_q_multi_lessthan1}
    p_{q}(x) = d[1+(q-1)(x-\mu_{q})^{\tr}\bm{B}(x-\mu_{q})]^{-\frac{1}{q-1}},
\end{equation}
and PDF is obtained by the following procedure as in the previous case as 
\begin{align}
p_{q}(x) &= \frac{1}{Z_{q}} \Big[1+\frac{q-1}{n+2-nq}(x-\mu_q)^{\tr}\Sigma_{q}^{-1}(x-\mu_q)\Big]^{-\frac{1}{q-1}}, \\ 
\text{with} \quad 
Z_{q} &= \Big[ \frac{(n+2)-nq}{q-1}\Big]^{\frac{n}{2}}|\Sigma_{q}|^{\frac{1}{2}}\pi^{\frac{n}{2}}\frac{\Gamma\big(1-
\frac{1}{q-1}\big)}{\Gamma \big(1-\frac{1}{q-1}+\frac{n}{2}\big)}, \ q < 1.
\end{align}

\noindent\textbf{Mean and Variance}\\
Here, we show mean and variance of the multivariate $q$-Gaussian in the case of $q>1$. For mean $\mu$, we have
\begin{align*}
\mu = 
\int x p_{q}(x)\mathrm{d}x &= d\int
x[1+b(q-1)(x-\mu_{q})^{\tr}\bm{B}(x-nu_{q})]^{-\frac{1}{q-1}} \mathrm{d}x \\
&= d \int z [1+(q-1)z^{\tr}\bm{B}z]^{-\frac{1}{q-1}} \mathrm{d}z + \mu_{q}
\end{align*}
To evaluate the first integral, we use the eigendecomposition of $\bm{B}$ and change of variable $t=\sqrt{q-1}\Lambda^{\frac{1}{2}}y)$, obtaining 
\begin{align*}
\int z [1+(q-1)z^{\tr}\bm{B}z]^{-\frac{1}{q-1}} \mathrm{d}z &= 
\sum_{i=1}^{n}\int y_{i}\bm{u}_{i}[1+(q-1)\sum_{k=1}^{n}\lambda_{k}y_{k}^{2}]^{-\frac{1}{q-1}}\mathrm{d}y \\
&= \sum_{i=1}^{n}\int [\lambda_{i}(q-1)]^{-\frac{1}{2}}\bm{u}_{i}t_{i}[1+\norm{t}^{2}]^{-\frac{1}{q-1}}\mathrm{d}t.
\end{align*}
From the symmetricity, we have 
\begin{align*}
\int t_{i} [1+\norm{t}^{2}]^{-\frac{1}{q-1}}\mathrm{d}t
&= t_{1} [1+\norm{t}^{2}]^{-\frac{1}{q-1}}\mathrm{d}t \\
& = \int r\cos\phi_{1} [1+r^{2}]^{-\frac{1}{q-1}} \mathrm{d}S_{n-1}(1)r^{n-1}\mathrm{d}r \\
&= \int [1+r^{2}]^{-\frac{1}{q-1}}r^{n}\mathrm{d}r \int \cos\phi_{1} \mathrm{d}S_{n-1}(1),
\end{align*}
where the second integral is converges to zero due to \eqref{eq:app_dS}. Thus, the integral above converges to zero if the first integral converges. To see this, we perform the change of variable $r^{2} = s$, which gives 
\begin{align*}
\int [1+r^{2}]^{-\frac{1}{q-1}}r^{n}\mathrm{d}r = \frac{1}{2}\int_{0}^{\infty}s^{\frac{n-1}{2}}(1+s)^{-\frac{1}{q-1}}\mathrm{d}s= B\Big(\frac{n}{2}+1, \frac{1}{q-1}-\frac{n+1}{2}\Big).
\end{align*}
This integral is convergent when $q<1+\frac{2}{n+1}$, which is tighter than \eqref{eq:multi_q_normal_cond}.
For the variance $\Sigma$, following the same procedure as we used to compute \eqref{eq:q_gauss_const_var}, 
\begin{align*}
\Sigma &= d 
\int(x-\mu_{q})(x-\mu_{q})^{\tr} p_{q}(x) \mathrm{d}x\\ 
&
= d(q-1)^{-\frac{n+2}{2}}|\bm{B}|^{-\frac{1}{2}}\sum_{i=1}^{n}\int \frac{\bm{u}_{i}\bm{u}_{i}^{\tr}}{\lambda_{i}}(1+\norm{t}^{2})^{-\frac{1}{q-1}}t_{i}^{2}\mathrm{d}t\\
&= d(q-1)^{-\frac{n+2}{2}}|\bm{B}|^{-\frac{1}{2}}\frac{A_{n-1}}{2n} B\Big(
\frac{1}{2}n+1, \frac{1}{q-1}-\frac{n+2}{2}\Big)\bm{B}^{-1}\\
&= \frac{n+2-nq}{n+4-(n+2)q}\Sigma_{q}, \ \text{with} \ q < 1 + \frac{2}{n+2},
\end{align*}
where the last inequality is required for convergence.
\begin{align}\label{eq:multi_q_statistics_range}
\mu = 
\begin{cases}
\mu_{q},  &1 <  q<1 + \frac{2}{n+1},\\
\text{undefined},  & 1 +\frac{2}{n+1} \leq q < 1 + \frac{2}{n},
\end{cases}
\quad 
\Sigma =
\begin{cases}
\frac{n+2-nq}{n+4-(n+2)q}\Sigma_{q}, \  & q < 1 + \frac{2}{n+2},\\
\infty, &1 + \frac{2}{n+2} \leq q < 1 +\frac{2}{n+1},\\
\text{undefined}, \quad & 1 +\frac{2}{n+1} \leq q < 1 + \frac{2}{n}.
\end{cases}
\end{align}
We note that substitution of $n=1$ recovers the the univariate case provided in \eqref{eq:summary_mu_sigma_uni}. 


\renewcommand{\theequation}{\thesubsection.\arabic{equation}}
\renewcommand{\thefigure}{\thesection.\arabic{figure}} 
\renewcommand{\thetable}{\thesection.\arabic{table}}
\setcounter{figure}{0} 
\setcounter{table}{0} 
\setcounter{equation}{0}
\subsection{Computing normalization constant of $\pi^{\ast}$}
This section provides a detailed derivation of how the normalization constant for the escort distribution of $\pi^{\ast}$ in section \ref{sec:q-gauss_policy}. From 
\eqref{eq:q_gauss_const_norm_q_multi}, \eqref{eq:multivariate_q_d-1} and  \eqref{eq:q_gaussB}, the normalization constant $C$ for $q$ escort distribution of multivariate $q$-Gaussian distribution is given as follows.  
\begin{align}
C = \Big[ \frac{(n+2)-nq}{q-1}\Big]^{\frac{n}{2}(1-q)}|\Sigma_{q}|^{\frac{1}{2}(1-q)}\pi^{\frac{n}{2}}\frac{\Gamma \big(\frac{q}{q-1}-\frac{n}{2}\big)}{\Gamma\big(
\frac{q}{q-1}\big)}\Bigg[\frac{\Gamma \big(\frac{1}{q-1}-\frac{n}{2}\Big)}{\Gamma\big(
\frac{1}{q-1}\big)}\Bigg]^{-q}.
\end{align}
From \eqref{eq:mu_q_sigma_q_q_Gausspolicy} we substitute $\Sigma_{q}$ and change $n$ to $n_{u}$, obtaining
\begin{align}
C = \frac{n_{u}+2-n_{u}q}{2} \Big( \frac{2[(q-1)\tilde{V}(x) + C\alpha]}{q-1}\Big)^{\frac{n_{u}}{2}(1-q)}|Q_{uu}^{-1}|^{\frac{1}{2}(1-q)}\pi^{\frac{n_{u}}{2}(1-q)} \Bigg[
\frac{\Gamma\big(\frac{1}{q-1}-\frac{n_{u}}{2})}{\Gamma\big(\frac{1}{q-1}\big)}
\Bigg]^{1-q}.
\end{align}
This equation recovers \eqref{eq:q_gauss_normalizer}.

\subsection{Existence of Escort Distribution}
Here, we examine the range of $q'$ introduced in the transformation in \eqref{eq:q_escort_trans_params}.
For the existence of finite mean and covariance, we have $1<q<1+2/(n+2)$. With this range of $q$, $q'=2-1/q$ takes the infimum
$$q' = 2-1/1 = 1$$
and supremum 
$$
2 - \frac{1}{1+\frac{2}{n+2}} = 2- \frac{n+2}{n+4} = 1 + \frac{2}{n+4} < 1 + \frac{2}{n+2}. 
$$
Therefore, when $q$ is in the appropriate range for the finite mean and covariance, the escort distribution also has finite mean and covariance.
\end{document}